\documentclass[10pt]{article}
\usepackage{amssymb, amsmath}
\usepackage{float,epsfig}
\usepackage{rotating}

\begin{document}

\newtheorem{theorem}{\bf Theorem}[section]
\newtheorem{proposition}[theorem]{\bf Proposition}
\newtheorem{definition}[theorem]{\bf Definition}
\newtheorem{corollary}[theorem]{\bf Corollary}
\newtheorem{example}[theorem]{\bf Example}
\newtheorem{exam}[theorem]{\bf Example}
\newtheorem{remark}[theorem]{\bf Remark}
\newtheorem{lemma}[theorem]{\bf Lemma}
\newcommand{\nrm}[1]{|\!|\!| {#1} |\!|\!|}

\newcommand{\ba}{\begin{array}}
\newcommand{\ea}{\end{array}}
\newcommand{\von}{\vskip 1ex}
\newcommand{\vone}{\vskip 2ex}
\newcommand{\vtwo}{\vskip 4ex}
\newcommand{\dm}[1]{ {\displaystyle{#1} } }

\newcommand{\be}{\begin{equation}}
\newcommand{\ee}{\end{equation}}
\newcommand{\beano}{\begin{eqnarray*}}
\newcommand{\eeano}{\end{eqnarray*}}
\newcommand{\inp}[2]{\langle {#1} ,\,{#2} \rangle}
\def\bmatrix#1{\left[ \begin{matrix} #1 \end{matrix} \right]}
\def \noin{\noindent}
\newcommand{\evenindex}{\Pi_e}

%\newcommand {\proof} {\par{\it Proof}. \ignorespaces}
%\newcommand {\eproof}
%      {\space
%        {\ \vbox{\hrule\hbox{\vrule height1.3ex\hskip0.8ex\vrule}\hrule}}
%        \par}

%%%%%%%%%%%%%%%%%%%%%%%%%%%%%%%%%%%%%%%%%%%%%%%%%%%%%%%%%%%%%%%%%%%%%%%%%%

\def \R{{\mathbb R}}
\def \C{{\mathbb C}}
\def \K{{\mathbb K}}
\def \J{{\mathbb J}}
\def \Lb{\mathrm{L}}

\def \T{{\mathbb T}}
\def \Pb{\mathrm{P}}
\def \N{{\mathbb N}}
\def \Ib{\mathrm{I}}
\def \Ls{{\Lambda}_{m-1}}
\def \Gb{\mathrm{G}}
\def \Hb{\mathrm{H}}
\def \Lam{{\Lambda_{m}}}
\def \Qb{\mathrm{Q}}
\def \Rb{\mathrm{R}}
\def \Mb{\mathrm{M}}
\def \norm{\nrm{\cdot}\equiv \nrm{\cdot}}

\def \P{{\mathbb P}_m(\C^{n\times n})}
\def \A{{{\mathbb P}_1(\C^{n\times n})}}
\def \H{{\mathbb H}}
\def \L{{\mathbb L}}
\def \G{{\mathcal G}}
\def \S{{\mathbb S}}
\def \sigmin{\sigma_{\min}}
\def \elam{\sigma_{\epsilon}}
\def \slam{\sigma^{\S}_{\epsilon}}
\def \Ib{\mathrm{I}}
\def \Tb{\mathrm{T}}
\def \d{{\delta}}

\def \Lb{\mathrm{L}}
\def \N{{\mathbb N}}
\def \Ls{{\Lambda}_{m-1}}
\def \Gb{\mathrm{G}}
\def \Hb{\mathrm{H}}
\def \Delta{\triangle}
\def \Rar{\Rightarrow}
\def \p{{\mathsf{p}(\lam; v)}}

\def \D{{\mathbb D}}

\def \tr{\mathrm{Tr}}
\def \cond{\mathrm{cond}}
\def \lam{\lambda}
\def \sig{\sigma}
\def \sign{\mathrm{sign}}

\def \ep{\epsilon}
\def \diag{\mathrm{diag}}
\def \rev{\mathrm{rev}}
\def \vec{\mathrm{vec}}

\def \sk{\mathsf{skew}}
\def \sy{\mathsf{sym}}
\def \en{\mathrm{even}}
\def \odd{\mathrm{odd}}
\def \rank{\mathrm{rank}}
\def \pf{{\bf Proof: }}
\def \dist{\mathrm{dist}}
\def \rar{\rightarrow}

\def \rank{\mathrm{rank}}
\def \pf{{\bf Proof: }}
\def \dist{\mathrm{dist}}
\def \Re{\mathsf{Re}}
\def \Im{\mathsf{Im}}
\def \re{\mathsf{re}}
\def \im{\mathsf{im}}

\def \sym{\mathsf{sym}}
\def \sksym{\mathsf{skew\mbox{-}sym}}
\def \odd{\mathrm{odd}}
\def \even{\mathrm{even}}
\def \herm{\mathsf{Herm}}
\def \skherm{\mathsf{skew\mbox{-}Herm}}
\def \str{\mathrm{ Struct}}
\def \eproof{$\blacksquare$}
\def \proof{\noin\pf}

\def \bS{{\bf S}}
\def \cA{{\cal A}}
\def \E{{\mathcal E}}
\def \X{{\mathcal X}}
\def \F{{\mathcal F}}
\def \tr{\mathrm{Tr}}
\def \range{\mathrm{Range}}

\def \pal{\mathrm{palindromic}}
\def \palpen{\mathrm{palindromic~~ pencil}}
\def \palpoly{\mathrm{palindromic~~ polynomial}}
\def \hodd{H\mbox{-}\odd}
\def \heven{H\mbox{-}\even}
%\def \herm{\mathrm{Hermitian}}
%\def \skherm{\mathrm{skew\mbox{-}Hermitian}}
%\def \str{\mathrm{ Struct}}

%%%%%%%%%%%%%%%%%%%%%%%%%%%%%%%%%%%%%%%%%%%%%%%%%%%%%%%%%%%%%%%%%%%%%%%%%%%%%%%%

\title{Backward errors and linearizations for palindromic matrix
polynomials}

\author{Bibhas Adhikari\thanks{Department of Mathematics, IIT Guwahati, Guwahati-781039,
 INDIA ({\tt bibhas.adhikari@gmail.com}).}}

\date{}

\maketitle

%\date

\vone \noin{\bf Abstract.}  We derive computable expressions of structured backward errors of approximate eigenelements of $*$-palindromic and $*$-anti-palindromic matrix polynomials. We also characterize minimal structured perturbations such that approximate eigenelements are exact eigenelements of the perturbed polynomials. We detect structure preserving linearizations which have almost no adverse effect on the structured backward errors of approximate eigenelements of the $*$-palindromic and $*$-anti-palindromic polynomials.

\vone \noin{\bf Keywords.} Structured backward error, palindromic matrix
polynomial, structured linearization

\vone\noin {\bf AMS subject classification(2000):} 
15A57, 65F35, 15A18.

\section{Introduction} A triple $(\lam,x,y)\in\C\times\C^n\times\C^n (x\neq 0, y\neq 0)$ is called an eigentriple of a polynomial $\Pb\in\P$ if \be\Pb(\lam)x=0 \,\, \mbox{and} \,\, y^H\Pb(\lam)=0, \ee where $\P$ denotes the space of matrix polynomials of the form $\Pb(z)=\sum_{j=0}^m z^jA_j, A_j\in\C^{n\times n}$ and $y^H$ is the conjugate transpose of $y.$ The nonzero vectors $x$ and $y$ are called the right and left eigenvectors of $\Pb$ corresponding to the eigenvalue $\lam,$ respectively. Given a polynomial $\Pb$ and a pair $(\lam,x)$ of $\Pb,$ the backward perturbation analysis deals with finding minimal perturbation $\Delta\Pb\in\P$ of $\Pb$ so that $(\lam,x)$ becomes an eigenpair of $\Pb+\Delta\Pb.$ If the coefficients of the given polynomial have certain distinctive structure sometimes it is necessary to find a minimal perturbation having the same structure as the original polynomial to preserve some properties (for example, eigensymmetry). 

In this paper we restrict our attention to regular matrix polynomials. We undertake a detailed backward perturbation analysis of $*$-palindromic and $*$-anti-palindromic matrix polynomials which we define in section \ref{sec1}. These polynomials arise in many applications such as in the study of rail traffic noise caused by high speed trains \cite{mackey2,mackey:thesis,schraoder:thesis,hillmm04}. Lately there has been a lot of interest generated into the development of structured preserving algorithms and the perturbation theory of palindromic polynomial eigenvalue problems \cite{ba1,chu:train,chu:palp,danchwat08,hillmm04,mmmm08,mackey2,schraoder:thesis}.

We denote the set of $*$-palindromic or $*$-anti-palindromic matrix polynomials by $\S\subset\P.$ We choose an appropriate norm $\nrm{\cdot}_M$ on $\P.$ Given a polynomial $\Pb \in \S$ and $(\lam, x) \in \C \times \C^n$ with $x^H x = 1$, we determine the structured backward error $\eta^\S_M (\lam, x, \Pb)$ of $(\lam, x)$ as an approximate right eigenpair of $\Pb \in \S$
and construct a polynomial $\Delta\Pb \in \S$ such that $\nrm{\Delta\Pb}_M = \eta^\S_M (\lam, x, \Pb)$ and $\Pb(\lam)x + \Delta\Pb(\lam)x = 0$. Moreover, we show that $\Delta\Pb$ is unique for the Frobenius norm on $\C^{n\times n}$ but there are infinitely many such $\Delta\Pb$ for the spectral norm on $\C^{n\times n}.$ Further, for the spectral norm, we show how to construct all such $\Delta\Pb$. A similar analysis undertaken in \cite{ba:poly} for certain other structures.

We mention that structured backward error multiplied with the structured condition number provides an approximate upper bound on the errors in the computed eigenelements. A detailed sensitivity analysis including explicit expression of structured condition number of eigenvalues of a variety of structured matrix polynomials including palindromic matrix polynomials has been investigated in \cite{ba3}. Thus structured backward errors derived in this paper will play an important role in the accuracy assessment of eigenelements of a $*$-palindromic/$*$-anti-palindromic matrix polynomial computed by structure preserving algorithms. 

Due to the lack of a genuine polynomial eigensolver, the common practice to solve a polynomial eigenvalue problem of degree $m$ is to solve an equivalent generalized eigenvalue problem of larger size. To be specific, an $n\times n$ polynomial $\Pb$ of degree $m$ is converted into an {\it equivalent} linear polynomial $\Lb(\lam)=\lam X + Y, \,\, X\in \C^{mn\times mn},Y\in \C^{mn\times mn}$ and a numerically backward stable algorithm is employed to compute the eigenelements of $\Lb.$ It is shown in \cite{mackey:thesis,mackey3} that a polynomial $\Pb\in\P$ can have infinitely many linearizations. In fact these linearizations form a vector space of dimension $m(m-1)n^2 + m.$ Analyzing backward error of approximate eigenpair and condition number of eigenvalues of a matrix polynomial Higham et al. \cite{higg1,higcon} have determined potential linearizations of a polynomial.

It is well known that $*$-palindromic/$*$-anti-palindromic matrix polynomials have certain eigensymmetry in the spectrum, and in the eigentriple as well \cite{mackey:thesis,mackey3,sbora}. Therefore to solve a palindromic polynomial eigenvalue problem it is very important to preserve those structures in the computed eigenelements. The structured linearizations which preserve the eigensymmetry of $*$-palindromic/$*$-anti-palindromic matrix polynomials have been constructed in \cite{mackey:thesis,mackey3}. Therefore computationally it is highly desirable to identify potential structured linearizations which are well-conditioned. By analyzing the structured condition number, a recipe of potential structured linearizations of a given $\Pb\in\S$ has been produced in \cite{ba3}.

With a view to analyzing accuracy of computed eigenelements of $*$-palindromic/$*$-anti-palindromic matrix polynomials we follow a similar procedure as developed in \cite{ba:poly} for a variety of structured polynomials including symmetric, skew-symmetric, even and odd. Indeed, we consider structured backward errors $\eta^\S_M(\lam, x, \Pb)$ of approximate eigenelements $(\lam,x)$ of $\Pb\in\S$ and structured backward errors $\eta^\S_M(\lam, \Lambda_{m-1}\otimes x,\Lb;v)$ of approximate eigenpair $(\lam, \Lambda_{m-1}\otimes x)$ of widely varying structured linearzations $\Lb$ of $\Pb,$ where $\Lambda_{m-1} := [ \lam^{m-1}, \ldots, \lam, 1]^T$ and $v$ is called the right ansatz vector, see \cite{mackey:thesis}. Further, we identify potential structured linearizations $\Lb$ of $\Pb$ for which $\eta^\S_M(\lam, \Lambda_{m-1}\otimes x, \Lb;v) \leq \alpha \eta^\S_M(\lam, x, \Pb),$ for some scalar $\alpha>0.$ Thus we identify structure preserving linearizations which have almost no adverse effect on the structured backward errors of approximate eigenelements of $*$-palindromic/$*$-anti-palindromic polynomials. We notice that the potential structured linearizations of $T$-palindromic matrix polynomials agree with those potential structured linearizations proposed in \cite{ba3} for $T$-palindromic matrix polynomials.

The rest of the paper is organized as follows. In section~\ref{sec1}, we review palindromic matrix polynomials and their spectral symmetries. In section~\ref{sec2}, we derive structured backward errors of approximate eigenpairs of palindromic matrix polynomials. In section~\ref{sec3},  we analyze  structured linearizations of palindromic matrix polynomials and identify potential structured linearizations.

\section{Eigensymmetry of palindromic matrix polynomials}\label{sec1}

%Any polynomial $\Pb\in\P$ is of the form
%$\Pb(z)=\sum_{j=0}^m z^j A_j.$ A pair $(\lam,x,y)$ is said
%to be an eigentriple of $\Pb,$ if
%
%\be\label{pal:eig}\Pb(\lam)x=0, \,\, y^H\Pb(\lam)=0.\ee
%
%An eigentriple $(\lam, x, y)$ is said to be normalized if $
%\|x\|_2 = 1 = \|y\|_2.$

A matrix polynomial $\Pb(z)=\sum_{j=0}^m z^j A_j\in\P$ is called 
$*$-palindromic or $*$-anti-palindromic if
\be\Pb^*(z)=z^m \Pb(1/z) \,\, \mbox{or} \,\, \Pb^*(z)=-z^m \Pb(1/z) \,\, \forall z\in\C\setminus\{0\} \ee respectively, where $\Pb^*(z)=\sum_{j=0}^m z^jA^*_j $ and $*\in\{T,H\}.$ Note that $A^T$ denotes the transpose of a matrix $A$ and the conjugate transpose of a matrix $A$ is denoted by $A^H.$ We denote the set of $*$-palindromic matrix polynomials by $\S_p$ and the set of $*$-anti-palindromic polynomials by $\S_{ap}.$ Unless otherwise stated we write $\S$ for both $\S_p$ and $\S_{ap}.$ Due to the structure of the
coefficients, the spectrum of a $*$-palindromic/$*$-anti-palindromic polynomial inherits a spectral symmetry. In fact if $\lam$ is an eigenvalue of $\Pb\in\S$ then $1/\lam^*$ is also an eigenvalue of $\Pb.$ This eigenvalue pairing $(\lam, 1/\lam^*)$ is known as the symplectic eigensymmetry. Table~\ref{table:paleigenpair} gives the eigensymmetry and structure of eigentriples of $*$-palindromic and $*$-anti-palindromic matrix polynomials.
\begin{table}[h]
\begin{center}\renewcommand{\arraystretch}{1.2}
\begin{tabular}{|c|c|c|}
  \hline
 % & \multicolumn{2}{|c|}{Polynomial $\Pb (\lambda)=\sum_{i=0}^m \lam^i A_i\in\S$}\\
% % & Pencil &  \\
%  \hline \hline
  \textbf{$\S$} & \textbf{eigenvalue pairing} & \textbf{eigentriple}  \\
  \hline \hline
    $T$-palindromic / $T$-antipalindromic & $(\lambda,1/\lam)$ & $(\lambda,x,\overline{y}),(1/\lam,y,\overline{x})$  \\
      % $T$-antipalindromic && \\
          \hline
    $H$-palindromic/$H$-anti-palindromic & $(\lambda,1/\overline{\lam})$ & $(\lambda,x,y),(1/\overline{\lam},y,x)$ \\
      %$H$-anti-palindromic  && \\
       \hline
\end{tabular}
\caption{\label{table:paleigenpair} Eigensymmetry of $*$-palindromic/$*$-anti-palindromic
polynomials.}
\end{center}
\end{table}

The results in Table~\ref{table:paleigenpair} follow from
\cite[Theorem 2.1]{ba1} by extending the arguments of matrix pencils to matrix polynomials. Note that if $\lam=0,$ that is, if $0$ is an eigenvalue of $\Pb\in\S$ then $\infty$ is an eigenvalue of $\Pb$ as well. In this paper we consider only finite eigenvalues, although an {\it infinity} eigenvalue can be analyzed by considering the {\it reversal} of the polynomial or by considering the homogeneous polynomial, see \cite{sa2,mackey:thesis,Tis00}. 

We now show that given $(\lam,x)\in\C\times\C^n$ with $\|x\|_2 =1$ and $\Pb\in\S,$ there always exists a polynomial $\Delta\Pb\in\S$ such that $(\Pb(\lam) + \Delta\Pb(\lam))x=0,$ that is, $(\lam,x)$ is an eigenpair of $\Pb+\Delta\Pb.$ For $x\in\C^n$ with $\|x\|_2=1,$ we define the projection $P_x :=I-xx^H.$ Throughout the paper, we follow the convention that $\Delta\Pb\in\P$ is of the form $\Delta\Pb(z)=\sum_{j=0}^m z^j\Delta A_j.$
\begin{theorem}\label{exist:tpal}
Let $\S\in\{\S_p,\S_{ap}\}$. Let $\Pb\in\S$ be given by $\Pb(z) = \sum_{j=0}^m z^j A_j.$ Suppose $(\lam,x)\in\C\times \C^n$ with $\|x\|_2=1.$ Set $r
:=-\Pb(\lam)x$ and $\Lam :=[1, \, \lam, \hdots, \, \lam^m]^T.$ Define
\beano \Delta A_j  &:=&\left\{%
\begin{array}{ll}
   -(x^{T}A_jx)\overline{x} x^{H}+\frac{1}{\|\Lam\|_2^2}\big[(\overline{\lambda})^{j}P_x^* rx^{H} +\ep (\overline{\lam})^{m-j} \, \overline{x} r^{T}P_x \big], & \hbox{if $* = T$} \\
   -(x^{H}A_jx)xx^{H}+\frac{1}{\|\Lam\|_2^2}\big[(\overline{\lambda})^{j}P_x rx^{H} + \ep\lam^{m-j}xr^{H}P_x\big], & \hbox{if $* = H$} \\
   \end{array}%
\right. \\
\Delta A_{m-j} &:=& \left\{%
\begin{array}{ll}
    \ep(\Delta A_j)^*, \, j=0:(m-1)/2 & \hbox{if $m$ is odd,} \\
    \ep(\Delta A_j)^*, \, j=0:m/2 & \hbox{if $m$ is even,} \\
\end{array}%
\right.\eeano where $\ep=1$ if $\S=\S_p,$ and $\ep=-1$ if $\S=\S_{ap}.$ Then $\Pb (\lam)x + \Delta \Pb(\lam)x =0$ and $\Delta\Pb\in\S.$ 
\end{theorem}
\noin \pf The proof is computational and is easy to
check.$\blacksquare$

\section{Structured backward error of approximate eigenpair}\label{sec2}

%Let $\Pb\in\S$ and $(\lam,x,y)$ be an approximate eigentriple of $\Pb.$ The structured backward %%%perturbation analysis mainly deals with finding the minimal perturbation $\Delta\Pb\in\S$ of $\Pb%%$ so that $(\lam,x,y)$ becomes an eigentriple of $\Pb+\Delta\Pb.$ In this paper we consider only %%the approximate right eigenpair $(\lam,x)\in\C\times \C^n,$ that is, $x$ is an approximate right %%eigenvector corresponding to an approximate eigenvalue $\lam.$ Unless otherwise stated throughout %the paper we always mean approximate right eigenpair $(\lam,x)$ by approximate eigenpair %%%%%%%%%%$(\lam,x).$ The similar arguments can be employed to derive the corresponding results for left %%%%eigenpair $(\lam,y).$

In this section we derive structured backward error of an approximate eigenpair $(\lam,x)$ of a polynomial $\Pb\in\S.$ The backward error of an approximate eigenpair $(\lam,x)$ is defined as the smallest perturbation by norm, $\Delta\Pb$ of $\Pb$ such that $(\lam,x)$ is an eigenpair of $\Pb+\Delta\Pb.$ Given $\Pb(z)=\sum_{j=0}^m z^jA_j$ we define norm in the following manner: \be\label{poly:norm}\nrm{\Pb}_M := \big(\sum_{j=0}^m \|A_j\|_M^2\big)^{1/2}, \,\, M\in\{F,2\}\ee where $M=F$ is the Frobenius norm and $M=2$ is the spectral norm. For a variety of norms on $\P$ see \cite{sa2}. 

Recall that a matrix polynomial $\Pb$ is called regular if $\det(\Pb(z)) \neq 0$ for some $z\in\C.$ Treating $(\lam, x)$ as an approximate eigenpair of a regular polynomial $\Pb$ we
define the backward error of $(\lam, x)$ by $$\eta_M(\lam,
x,\Pb) := \min\limits_{\Delta \Pb\in\P} \big\{ \nrm{\Delta\Pb}_M :
 \Pb(\lam)x +\Delta\Pb(\lam)x = 0\big\}, \,\, M\in\{F,2\}.$$
Setting $r :=-\Pb(\lam)x,$ we have the explicit formula \cite{ba:poly} \be\label{eq:unbckerr}\eta_M(\lam, x, \Pb)   =\|r\|_2/\|x\|_2\|\Lam\|_2.\ee An explicit formula of backward error is obtained by Tisseur \cite{Tis00} for a different class of norms on $\P.$ See also \cite{sbora}.

Next assume that $\Pb\in\S$ is a regular polynomial. Then we define the structured
backward error of an approximate eigenpair $(\lam, x)$ by \be\label{def:strerr}\eta_M^\S(\lam,x,\Pb) := \min\limits_{\Delta \Pb\in\S} \big\{ \nrm{\Delta\Pb}_M :
 \Pb(\lam)x +\Delta\Pb(\lam)x = 0\big\}, \,\, M\in\{F,2\}.\ee By Theorem~\ref{exist:tpal} it is easy to see that $\eta_M^\S(\lam,x,\Pb)< \infty$ and $\eta_M(\lam,x,\Pb) \leq \eta_M^\S(\lam,x,\Pb).$

To derive $\eta_M^\S(\lam,x,\Pb)$ corresponding to $M=2$ we use Davis-Kahan-Weinberger norm-preserving dilation theorem (DKW in short) which we state below.  
\begin{theorem}[Davis-Kahan-Weinberger, \cite{D}] \label{dkw}
Let $A,B,C$ and $D$ are matrices of appropriate sizes. Let $A, B, C$ satisfy $ \left\|\bmatrix{ A\\
B}\right\|_2 = \mu$ and $ \left\|\bmatrix{ A & C}\right\|_2
= \mu.$ Then   there exists $D$ such that $\left\|\bmatrix{  A & C \\
  B & D}\right\|_2 = \mu.$ All $D$ which have this property are
exactly those of the form \be\label{dkw:d} D= - KA^HL + \mu(I-KK^H)^{1/2}Z(I-L^HL)^{1/2},\ee where $K^H :=(\mu^2 \Ib -A^HA )^{-1/2}B^H,~ L :=(\mu^2 \Ib -AA^H)^{-1/2}C$ and $Z$ is an
arbitrary contraction, that is, $\|Z\|_2 \leq 1.$ 
\end{theorem} For a more general version of DKW Theorem see \cite{D}. We use DKW Theorem in the subsequent development by setting $Z=0$ in (\ref{dkw:d}) to avoid cumbersome calculations.

Let $\Lam :=[1, \, \lam, \, \ldots \, \lam^m]^T, \lam\in\C.$ To determine structured backward error in a convenient manner we use the projection operators $\Pi_s $ which were introduced in \cite{ba3} to determine the structured condition number of eigenvalues of $*$-palindromic/$*$-anti-palindromic matrix polynomials. $\Pi_s, s\in\{+,-\}$ is defined by
\begin{equation} \label{eq:symmetricprojector}
    \Pi_{\pm}(\Lam) :=
 \left\{\!\!\!\!
 \begin{array}{ll}
     \big[ \frac{\lambda^m \pm 1}{\sqrt{2}},
%\,\frac{ \lambda^{m-1} +  \lambda}{\sqrt{2}},
\,\ldots,\, \frac{\lambda^{m/2+1} \pm
\lambda^{m/2-1}}{\sqrt{2}},\,  \frac{\lambda^{m/2}\pm
\lambda^{m/2}}{2} \big]^T  & \text{if $m$ is even}, \\{}
          \big[ \frac{\lambda^m \pm 1}{\sqrt{2}},
\,\ldots,\, \frac{\lambda^{(m+1)/2} \pm \lambda^{(m-1)/2}}{\sqrt{2}} \big]^T, & \text{if $m$ is odd.} \\
 \end{array} \right.
\end{equation} Now we state some basic properties of $\Pi_+$ and $\Pi_-$ that will be used in the subsequent development. It is straightforward to check that the following relations hold.
\begin{itemize} \item $\|\Pi_+(\Lam)\|_2^2 - \|\Pi_-(\Lam)\|_2^2 =  \left\{\!\!\!\!
 \begin{array}{ll}
 \sum_{j=0}^{(m-2)/2} 2\re((\overline{\lam})^j\lam^{m-j}) + |\lam^{m/2}|^2, & \text{if $m$ is even} \\{}
\sum_{j=0}^{(m-1)/2} 2\re((\overline{\lam})^j\lam^{m-j}), & \text{if $m$ is odd.}\\
\end{array} \right.$ \item $\|\Pi_+(\Lam)\|_2^2 + \|\Pi_-(\Lam)\|_2^2 = \|\Lam\|_2^2$. \item $2\|\Pi_+(\Lam) + \Pi_-(\Lam)\|_2^2 =  \left\{\!\!\!\!
 \begin{array}{ll}
 \sum_{j=0}^{m/2} |\lam^{m-j}|^2, & \text{if $m$ is even} \\{}
\sum_{j=0}^{(m-1)/2} |\lam^{m-j}|^2, & \text{if $m$ is odd.}\\
\end{array} \right.$ \item $2\|\Pi_+(\Lam) - \Pi_-(\Lam)\|_2^2 =  \left\{\!\!\!\!
 \begin{array}{ll}
 \sum_{j=0}^{m/2} |\lam^{j}|^2, & \text{if $m$ is even} \\{}
\sum_{j=0}^{(m-1)/2} |\lam^{j}|^2, & \text{if $m$ is odd.}\\
\end{array} \right.$ \end{itemize} 

\subsection{$T$-palindromic and $T$-anti-palindromic matrix polynomials}

Now we derive structured backward error of approximate eigenpair of $T$-palindromic and $T$-anti-palindromic matrix polynomials. Recall that a polynomial $\Pb(z)=\sum_{j=0}^m z^j A_j$ is $T$-palindromic if $A_j^T= A_{m-j},$ and $T$-anti-palindromic if $A_j^T=- A_{m-j}.$ The set of $T$-palindromic and $T$-anti-palindromic matrix polynomials is denoted by $\S_p$ and $\S_{ap}$ respectively.

\begin{theorem}\label{bck:tpal}
Let $\S\in\{\S_p, \S_{ap}\}.$ Let $\Pb\in\S$ be given by $\Pb(z)=\sum_{j=0}^m z^jA_j.$ Let $(\lam,x)$ with $\|x\|_2=1$ be an approximate eigenpair of $\Pb.$ Set $r :=-\Pb(\lam)x.$ Then we have $$\eta^\S_M(\lam,x,\Pb) = \big( a_M^{\S}(\lam)\|r\|_2^2 + b_M^{\S}(\lam)|x^Tr|^2 \big)^{1/2}, \,\, M\in\{F,2\}$$ where $a_M^\S(\lam)$ and $b_M^\S(\lam)$ are given by 
%\begin{table}[hs]
\begin{center}\renewcommand{\arraystretch}{1.2}
\begin{tabular}{|l|c|c|c|c|}
  \hline
   $m$ & $\lam$ & $\S$ & $a_M^\S(\lam), \,\, M=F$ & $b_M^\S(\lam), \,\, M=F$\\
\hline \hline
   %after \\: \hline or \cline{col1-col2} \cline{col3-col4} ...
 odd & $\lam\neq -1$ & $\S_{p}$ &$\frac{2}{\|\Lam\|_2^2}$ & $\frac{1}{\|\Pi_{s}(\Lam)\|_2^2}-\frac{2}{\|\Lam\|_2^2}$ \\\cline{2-3}
& $\lam\neq 1$ & $\S_{ap}$ && \\\cline{2-5}
 %\hline\cline{2-3}
 & $\lam=-1$ & $\S_{p}$ & $\frac{2}{\|\Lam\|_2^2}$ & $0$ \\ \cline{2-3}
& $\lam=1$ & $\S_{ap}$ && \\
\hline
even & $\lam = \pm 1$ & $\S_{p}$ & $\frac{2}{\|\Lam\|_2^2}$ & $-\frac{1}{\|\Lam\|_2^2}$ \\\cline{2-5}
%\hline
 & $\lam\neq \pm 1$ & $\S_{p}$ & $\frac{2}{\|\Lam\|_2^2}$ & $\frac{4\|\Pi_s(\Lam)\|_2^2 - 3\ep|\lam^{m/2}|^2}{\big(2\|\Pi_s(\Lam)\|_2^2 - \ep|\lam^{m/2}|^2\big)^2} - \frac{2}{\|\Lam\|_2^2}$ \\\cline{2-3}
& $\lam\neq 1$ & $\S_{ap}$ && \\
\hline 
\end{tabular}
\end{center}
\begin{center}\renewcommand{\arraystretch}{1.2}
\begin{tabular}{|l|c|c|c|c|}
\hline
$m$ & $\lam$ & $\S$ & $a_M^\S(\lam), \,\, M=2$ & $b_M^\S(\lam), \,\, M=2$\\
\hline \hline
 odd & $|\lam|>1$ & $\S$ & $\frac{4\|\Pi_+(\Lam) + \Pi_-(\Lam)\|_2^2 }{\|\Lam\|_2^4}$ & $\frac{\|\Lam\|_2^4 -4\|\Pi_+(\Lam) + \Pi_-(\Lam)\|_2^2 \|\Pi_s(\Lam)\|_2^2}{\|\Lam\|_2^4 \|\Pi_s(\Lam)\|_2^2}$\\ 
%&& $\S_{ap}$ && \\\cline{2-5}
 & $|\lam|\leq 1$ & $\S$  & $\frac{4\|\Pi_+(\Lam) - \Pi_-(\Lam)\|_2^2 }{\|\Lam\|_2^4}$ & $\frac{\|\Lam\|_2^4 -4\|\Pi_+(\Lam) - \Pi_-(\Lam)\|_2^2 \|\Pi_s(\Lam)\|_2^2}{\|\Lam\|_2^4 \|\Pi_s(\Lam)\|_2^2}$\\
%&& $\S_{ap}$ && \\
 \hline
even & $\lam = \pm 1$ & $\S_{p}$ & $\frac{1}{\|\Lam\|_2^2}$ & $0$ \\\cline{2-3}
& $\lam=1$ & $\S_{ap}$ && \\\cline{2-5}
& $|\lam|>1$ & $\S$ &$\frac{4\|\Pi_+(\Lam) + \Pi_-(\Lam)\|_2^2- |\lam^{m/2}|^2 }{\|\Lam\|_2^4 }$  & $\big[\frac{4\|\Pi_s(\Lam)\|_2^2 - 3\ep|\lam^{m/2}|^2}{\big(2\|\Pi_s(\Lam)\|_2^2 - \ep|\lam^{m/2}|^2\big)^2} - $ \\
&&  && $\frac{4\|\Pi_+(\Lam) + \Pi_-(\Lam)\|_2^2 - |\lam^{m/2}|^2}{\|\Lam\|_2^4}\big]$ \\\cline{2-5}
%&& $\S_{ap}$ && \\\cline{2-5}
& $|\lam|\leq 1$ & $\S$ &$\frac{4\|\Pi_+(\Lam) - \Pi_-(\Lam)\|_2^2 - |\lam^{m/2}|^2}{\|\Lam\|_2^4}  $ & $\big[\frac{4\|\Pi_s(\Lam)\|_2^2 - 3\ep|\lam^{m/2}|^2}{\big(2\|\Pi_s(\Lam)\|_2^2 - \ep|\lam^{m/2}|^2\big)^2} - $ \\
&& && $\frac{4\|\Pi_+(\Lam) - \Pi_-(\Lam)\|_2^2 - |\lam^{m/2}|^2}{\|\Lam\|_2^4}\big]$ \\
%&& $\S_{ap}$ && \\
\hline
\end{tabular}
%\caption{\label{table:tpalbckerr} }
\end{center}
%\end{table} 
where $s=+,\ep=1$ if $\S=\S_p,$ and $s=-,\ep=0$ if $\S=\S_{ap}.$
%\begin{table}[ht]
%\begin{center}\renewcommand{\arraystretch}{1.0}
%\begin{tabular}{|l|l|l|l|}
 % \hline
%   \multicolumn{2}{|c|}{$a_M^\S(\lam)$} & \multicolumn{2}{|c|}{$b_M^\S(\lam)$}\\
%\hline \hline
  %$M=F$ & $M=2$ & $M=F$ & $M=2$ \\
 % \hline \hline
    %after \\: \hline or \cline{col1-col2} \cline{col3-col4} ...
 %$\frac{2}{\|\Lam\|_2^2}$  & $\frac{4\|\Pi_+(\Lam) + \Pi_-(\Lam)\|_2^2 }{\|\Lam\|_2^4}$ & $\frac%%%{\|\Pi_{-}(\Lam)\|_2^2 - \|\Pi_{+}(\Lam)\|_2^2}{\|\Lam\|_2^2 \|\Pi_{+}(\Lam)\|_2^2}$ & $\frac{\|%%%\Lam\|_2^4 -4\|\Pi_+(\Lam) + \Pi_-(\Lam)\|_2^2 \|\Pi_+(\Lam)\|_2^2}{\|\Lam\|_2^4 \|\Pi_+(\Lam)\|%%%_2^2}$ \\
 %if $\lam\neq -1$ & if $|\lam|>1$ & if $\lam\neq -1$ & if $|\lam|> 1$ \\\cline{1-4}
 % $\frac{2}{(m+1)}$ & $\frac{4\|\Pi_+(\Lam) - \Pi_-(\Lam)\|_2^2 }{\|\Lam\|_2^4}$ & $0$& $\frac{\|%%%\Lam\|_2^4 -4\|\Pi_+(\Lam) - \Pi_-(\Lam)\|_2^2 \|\Pi_+(\Lam)\|_2^2}{\|\Lam\|_2^4 \|\Pi_+(\Lam)\|%%%%_2^2}$ \\
 %if $\lam=-1$ & if $|\lam|\leq 1$ & if $\lam=-1$ & if $|\lam|\leq 1$\\    
%\hline
%\end{tabular}
%\caption{\label{table:tpalbckerr} }
%\end{center}
%\end{table} 
\end{theorem}

\noin\pf First suppose that $m$ is even. Then by Theorem~\ref{exist:tpal}, note that there aways exists a polynomial $\Delta \Pb\in\S$ which satisfies $\Delta \Pb(\lam)x + \Pb(\lam)x=0.$ Consequently $\eta_M^{\S} (\lam,x,\Pb) < \infty.$ Let $Q=[x, \,\,Q_1]$ be a unitary matrix, where $x\in\C^n$ is given and $Q_1\in\C^{n\times (n-1)}$ is an isometry so that $Q_1^Hx=0.$ 

Let $\S=\S_p.$ Define \be\label{eq1tpal}\Delta A_j := \overline{Q}\bmatrix{a_{jj} & a_j^T \\
  b_j & X_j} Q^H,~ \Delta A_{m/2} := \overline{Q}\bmatrix{a_{(m/2)(m/2)} & a_{m/2}^T \\
  a_{m/2} & X_{m/2}}Q^H, j=0:(m-2)/2\ee
and $\Delta A_{m-j} = (\Delta A_j)^T.$ Now $\Delta \Pb(\lam)x + \Pb(\lam)x=0$ yields $\Delta\Pb(\lam)x=-\Pb(\lam)x=r$(say). Therefore by (\ref{eq1tpal}) we have  
$$\bmatrix{\sum_{j=0}^{m} \lam^j a_{jj} \\ \sum_{j=0}^{(m-2)/2} \lam^j b_j +\sum_{j=0}^{(m-2)/2} \lam^{m-j} a_j +\lam^{m/2} a_{m/2}}=\bmatrix{x^Tr \\ Q_1^Tr }.$$ By Lemma \ref{app:2} the minimum norm solution of $\sum_{j=0}^{(m-2)/2} \lam^j b_j +\sum_{j=0}^{(m-2)/2} \lam^{m-j}
a_j +\lam^{m/2} a_{m/2}=Q_1^Tr$ is given by $$ b_j =
\frac{(\overline{\lam})^j}{\|\Lam\|_2^2}Q_1^Tr,\,
    a_j = \frac{(\overline{\lam})^{m-j}}{\|\Lam\|_2^2}Q_1^Tr,\,j=0:(m-2)/2, \,\,
    a_{m/2} = \frac{(\overline{\lam})^{m/2}}{\|\Lam\|_2^2}Q_1^Tr.$$ Further, by Lemma \ref{app:1}, \ref{it:4} the minimum norm solution of $\sum_{j=0}^{m}\lam^j a_{jj}=x^Tr$ is given by $$a_{jj} =\left\{%
\begin{array}{ll}
\frac{\lam^j}{\|\Lam\|_2^2}x^{T}r & \hbox{if $\lam=\pm 1,$}\\
\frac{(\overline{\lam})^j +(\overline{\lam})^{m-j}}{2\|\Pi_+(\Lam)\|_2^2 - |\lam^{m/2}|^2}x^Tr & \hbox{if $\lam\neq \pm 1,$}\\
\end{array}%
\right.  a_{(m/2)(m/2)} =\left\{%
\begin{array}{ll}
\frac{\lam^{m/2}}{\|\Lam\|_2^2}x^{T}r \hfill{\mbox{if}\, \lam=\pm 1,}\\
\frac{(\overline{\lam})^{m/2}x^Tr}{2\|\Pi_+(\Lam)\|_2^2 - |\lam^{m/2}|^2} \hfill{\mbox{if}\, \lam\neq \pm 1,}\\
\end{array}%
\right. $$ where $j=0:(m-2)/2.$ Thus we have \be\label{apm1f}\Delta A_j = \overline{Q}\bmatrix{
  \frac{\lam^j}{\|\Lam\|_2^2}x^{T}r & \frac{\lam^{m-j}}{\|\Lam\|_2^2}(Q^{T}_{1}r)^{T}\\
  \frac{\lam^{j}}{\|\Lam\|_2^2}Q^{T}_{1}r & X_j }Q^{H}, j=0:m/2, \Delta A_{m-j} = (\Delta A_j)^T\ee whenever $\lam=\pm 1,$ and \be\label{apmn1f} \Delta A_j = \overline{Q}\bmatrix{
  \frac{(\overline{\lam})^j +(\overline{\lam})^{m-j}}{2\|\Pi_+(\Lam)\|_2^2-|\lam^{m/2}|^2}x^Tr &
 \frac{(\overline{\lam})^{m-j}}{\|\Lam\|_2^2}(Q^{T}_{1}r)^{T} \\
  \frac{(\overline{\lam})^{j}}{\|\Lam\|_2^2}Q^{T}_{1}r & X_j }Q^{H}, \, j=0:(m-2)/2, \ee\be\label{apmn1mf} \Delta A_{m/2} =  \overline{Q}\bmatrix{
  \frac{(\overline{\lam})^{m/2}}{2\|\Pi_+(\Lam)\|_2^2-|\lam^{m/2}|^2}x^Tr &
  \frac{(\overline{\lam})^{m/2}}{\|\Lam\|_2^2}(Q^{T}_{1}r)^{T} \\
  \frac{(\overline{\lam})^{m/2}}{\|\Lam\|_2^2}Q^{T}_{1}r & X_{m/2} } Q^{H},\, \Delta A_{m-j} = (\Delta A_j)^T\ee whenever $\lam\neq \pm 1.$ Setting $X_j=0, j=0:m,$ and using the fact that 
\be\label{qtr} \|Q^T_1r\|_2^2 = \|\overline{Q}Q^T_1r\|_2^2 =\|(I-\overline{x}x^T)r\|_2^2=\|r\|_2^2-|x^Tr|^2, \ee we obtain 
$$\eta_F^{\S_p}(\lam,x,\Pb)=\left\{%
\begin{array}{ll}
\frac{1}{\sqrt{m+1}}~\sqrt{2\|r\|_{2}^{2}-~|x^{T}r|^{2}} & \hfill{\mbox{if}\, \lam=\pm 1} \\
\sqrt{\frac{2}{\|\Lam\|_2^2}\|r\|_2^2 + \big(\frac{4\|\Pi_+(\Lam)\|_2^2 - 3|\lam^{m/2}|^2}{\big(2\|\Pi_+(\Lam)\|_2^2 - |\lam^{m/2}|^2\big)^2} - \frac{2}{\|\Lam\|_2^2}\big)|x^Tr|^2}& \hfill{\mbox{if} \, \lam\neq \pm 1.}\\
\end{array}%
\right.$$

Now we derive $\eta_2^{\S_p}(\lam,x,\Pb),$ by using DKW Theorem. If $\lam=\pm 1,$ by (\ref{apm1f}) and Theorem~\ref{dkw} we have $ \mu_{\Delta A_j}=\|\Delta A_j\|_2 =\frac{\|r\|_2}{\|\Lam\|_2^4},~j=0:m,$ given by \be\label{x:pm1} X_j=-\frac{\lam^{m-j} \overline{x^Tr}~Q_1^Tr(Q_1^Tr)^T}{\|\Lam\|_2^2~(\|r\|_2^2-|x^Tr|^2)}, \, X_{m/2}=-\frac{\lam^{m/2}\overline{x^Tr}~Q_1^Tr(Q_1^Tr)^T}{\|\Lam\|_2^2~(\|r\|_2^2-|x^Tr|^2)}\ee where $~j=0:(m-2)/2,$ and $X_{m-j}=X_j^T.$ If $\lam\neq \pm 1,$ by (\ref{apmn1f}), (\ref{apmn1mf}) and Theorem~\ref{dkw} we have \beano \mu_{\Delta A_j} &=& \|\Delta A_j\|_2 =\left\{%
\begin{array}{ll}
    \sqrt{\frac{|\lam^j + \lam^{m-j}|^2~|x^Tr|^2}{\big(2\|\Pi_+(\Lam)\|_2^2-|\lam^{m/2}|^2\big)^2} + \frac{|\lam^{m-j}|^2~(\|r\|_2^2-|x^Tr|^2)}{\|\Lam\|_2^4}}, & \hbox{if $|\lam|>1$} \\[8pt]
   \sqrt{\frac{|\lam^j + \lam^{m-j}|^2~|x^Tr|^2}{\big(2\|\Pi_+(\Lam)\|_2^2-|\lam^{m/2}|^2\big)^2}+ \frac{|\lam^{j}|^2~(\|r\|_2^2-|x^Tr|^2)}{\|\Lam\|_2^4}}, & \hbox{if $|\lam|\leq 1$} \\
\end{array}%
\right. \\ \mu_{\Delta A_{m/2}} &=& \|\Delta A_{m/2}\|_2 =\sqrt{\frac{|\lam^{m/2}
|^2~|x^Tr|^2}{\big(2\|\Pi_+(\Lam)\|_2^2-|\lam^{m/2}|^2\big)^2}+
\frac{|\lam^{m/2}|^2~(\|r\|_2^2-|x^Tr|^2)}{\|\Lam\|_2^4}}\eeano given by
\be X_j = \left\{%
\begin{array}{ll}
    -\frac{(|\lam^j|^2(\overline{\lam})^{m-j}+|\lam^{m-j}|^2(\overline{\lam})^j~\overline{x^Tr}~Q_1^Tr
(Q_1^Tr)^T}{\big(2\|\Pi_+(\Lam)\|_2^2-|\lam^{m/2}|^2\big)~|\lam^{m-j}|^2~(\|r\|_2^2-|x^Tr|^2)}, & \hbox{if $|\lam|>1$,} \\[10pt]
   -\frac{(|\lam^j|^2(\overline{\lam})^{m-j}+|\lam^{2m-j}|^2(\overline{\lam})^j~\overline{x^Tr}~Q_1^Tr
(Q_1^Tr)^T}{\big(2\|\Pi_+(\Lam)\|_2^2-|\lam^{m/2}|^2\big)~|\lam^{j}|^2~(\|r\|_2^2-|x^Tr|^2)}, & \hbox{if $|\lam|\leq 1$,} \\
\end{array}%
\right.\ee \be\label{x:npm11}
 X_{m/2} = -\frac{(\overline{\lam})^{m/2}~\overline{x^Tr}~Q_1^Tr(Q_1^Tr)^T}{\big(2\|\Pi_+(\Lam)\|_2^2-|\lam^{m/2}|^2\big)~(\|r\|_2^2-|x^Tr|^2)},\ee where $j=0:(m-2)/2.$ Consequently we have $\eta_2^{\S_p}(\lam,x,\Pb)=\frac{\|r\|_{2}}{\sqrt{m+1}}$ if $\lam=\pm 1,$ and $$\big( \eta_2^{\S_p}(\lam,x,\Pb) \big)^2 = \left\{%
\begin{array}{ll}
\frac{4\|\Pi_+(\Lam) + \Pi_-(\Lam)\|_2^2- |\lam^{m/2}|^2 }{\|\Lam\|_2^4 }\|r\|_2^2 + \big[\frac{4\|\Pi_+(\Lam)\|_2^2 - 3|\lam^{m/2}|^2}{\big(2\|\Pi_+(\Lam)\|_2^2 - |\lam^{m/2}|^2\big)^2} \\\hfill{ - \frac{4\|\Pi_+(\Lam) + \Pi_-(\Lam)\|_2^2 - |\lam^{m/2}|^2}{\|\Lam\|_2^4}\big]|x^Tr|^2 \, \mbox{if} \, |\lam|>1,} \\
\frac{4\|\Pi_+(\Lam) - \Pi_-(\Lam)\|_2^2 - |\lam^{m/2}|^2}{\|\Lam\|_2^4} \|r\|_2^2 +\big[ \frac{4\|\Pi_+(\Lam)\|_2^2 - 3|\lam^{m/2}|^2}{\big(2\|\Pi_+(\Lam)\|_2^2 - |\lam^{m/2}|^2\big)^2}\\ \hfill{ - \frac{4\|\Pi_+(\Lam) - \Pi_-(\Lam)\|_2^2 - |\lam^{m/2}|^2}{\|\Lam\|_2^4}\big]|x^Tr|^2 \, \mbox{if} \, |\lam|\leq 1, \lam\neq \pm 1.}
\end{array}%
\right.$$ Note that if $|x^Tr|=\|r\|_2,$ then $\|Q_1^Tr\|_2=0$. In such a case,
considering $X_j=0$ we obtain the desired results.  

Next let $\S=\S_{ap}.$ Define \be\label{eq1tantipal}\Delta A_j := \overline{Q}\bmatrix{a_{jj} & a_j^T \\b_j & X_j} Q^H,~ \Delta A_{m/2} := \overline{Q}\bmatrix{0 & - a_{m/2}^T \\
  a_{m/2} & X_{m/2}}Q^H, j=0:(m-2)/2\ee
and $\Delta A_{m-j} = -(\Delta A_j)^T$ and $Q=[x, \,\, Q_1]$ is a unitary matrix. Therefore we have $$\bmatrix{
  \sum_{j=0}^{(m-2)/2} \lam^j a_{jj} - \sum_{j=0}^{(m-2)/2} \lam^{m-j} a_{jj}\\
  \sum_{j=0}^{(m-2)/2} \lam^j b_j + a_{m/2} \lam^{m/2} - \sum_{j=0}^{(m-2)/2} \lam^{m-j} a_j}= \bmatrix{ x^Tr \\   Q_1^Tr }.$$ Note that $a_{jj}=0$ whenever $j=m/2,$ since $(A_{m/2})^T=-A_{m/2}.$ By Lemma \ref{app:2}, the minimum norm solution of
$\sum_{j=0}^{(m-2)/2} \lam^j b_j + a_{m/2} \lam^{m/2} -
\sum_{j=0}^{(m-2)/2} \lam^{m-j} a_j = Q_1^T r$ is given by $$
    b_j = \frac{(\overline{\lam})^j Q_1^Tr}{\|\Lam\|_2^2}, \,\,
    a_j = -\frac{(\overline{\lam})^{m-j}Q_1^Tr}{\|\Lam\|_2^2},
    \,\, a_{m/2} = \frac{(\overline{\lam})^{m/2} Q_1^T r}{\|\Lam\|_2^2}, j=0:(m-2)/2.$$ Also note that $x^Tr=0$ if $\lam=1.$ Therefore by Lemma \ref{app:1}, \ref{it:5}, the minimum norm solution of $\sum_{j=0}^{(m-2)/2} \lam^j a_{jj} - \sum_{j=0}^{(m-2)/2} \lam^{m-j} a_{jj} = x^Tr$ is given by $$a_{jj} = \left\{%
\begin{array}{ll}
 0 \hfill{\mbox{if} \, \lam=1}\\
 \frac{(\overline{\lam})^j - (\overline{\lam})^{m-j}}{
2\|\Pi_-(\Lam) \|_2^2}x^Tr \,  \hfill{\mbox{if} \, \lam\neq 1}\\
\end{array}%
\right.$$ where $j=(m-2)/2.$ Therefore by (\ref{eq1tantipal}) we have \be \Delta A_j = \left\{%
\begin{array}{ll}
\overline{Q}\bmatrix{ 0 & -\frac{(Q_1^Tr)^T}{\|\Lam\|_2^2} \\ \frac{Q_1^Tr}{\|\Lam\|_2^2} & X_j}Q^H \hfill{\mbox{if} \, \lam=1}\\
\overline{Q}\bmatrix{\frac{(\overline{\lam})^j -
(\overline{\lam})^{m-j}}{ 2\| \Pi_-(\Lam) \|_2^2}x^Tr & -\frac{(\overline{\lam})^{m-j}(Q_1^Tr)^T}{\|\Lam\|_2^2}\\ \frac{(\overline{\lam})^j Q_1^Tr}{\|\Lam\|_2^2} & X_j }Q^H \hfill{\mbox{if} \, \lam\neq 1}
\end{array}%
\right. \ee \be \Delta A_{m/2} =\overline{Q} \bmatrix{0 & -\frac{(\overline{\lam})^{m/2}(Q_1^Tr)^T}{\|\Lam\|_2^2} \\ \frac{(\overline{\lam})^{m/2} Q_1^Tr}{\|\Lam\|_2^2} & X_{m/2}}Q^H \, \mbox{if} \, \lam\neq 1,\ee where $j=0:(m-2)/2$ and $\Delta A_{m-j}=(\Delta A_j)^T.$ Now setting $X_j=0,$ we obtain the desired result for $M=F.$ Further, by employing DKW Theorem \ref{dkw} and following a similar arguments as that in the case of $\S=\S_p,$ we obtain $\eta^{\S_{ap}}_2(\lam,x,\Pb).$

Next consider $m$ be odd. Let $\S=\S_p.$ Define $$ \widetilde{\Delta A}_j :=Q^T\Delta AQ=\bmatrix{
  a_{jj} & a_j^T \\
  b_j & X_j }~\mbox{and} ~~(\Delta A_j)^T=\Delta A_{m-j}, j =0:(m-1)/2,$$ where $Q=[x, \, Q_1], Q_1\in \C^{n \times (n-1)}$ is a unitary matrix defined as above. Consequently we have $$\bmatrix{
  \sum_{j=0}^{m} \lam^j a_{jj} \\
  \sum_{j=0}^{(m-1)/2} \lam^j b_j +\sum_{j=0}^{(m-1)/2} \lam^{m-j} a_j }=\bmatrix{
  x^Tr \\
  Q_1^Tr }.$$ Note that $x^Tr=0$ if $\lam=-1.$ Then by Lemma \ref{app:2} and Lemma \ref{app:1}, \ref{it:3} the minimum norm solutions of $\sum_{j=0}^{(m-1)/2} \lam^j b_j + 
\sum_{j=0}^{(m-1)/2} \lam^{m-j} a_j =Q_1^Tr$ and $\sum_{j=0}^{m} \lam^j a_{jj}=x^Tr$ are given by 
$$b_j=\frac{(\overline{\lam})^j}{\|\Lam\|_2^2} Q_1^T r, a_j= \frac{(\overline{\lam})^{m-j}}{\|\Lam\|_2^2} Q_1^Tr, a_{jj} = \left\{%
\begin{array}{ll}
0 \hfill{\, \mbox{if} \, \lam=-1,}\\
\frac{(\overline{\lam})^{j}+(\overline{\lam})^{m-j}}{2\|\Pi_+(\Lam)\|_2^2}x^Tr \hfill{\, \mbox{if} \, \lam\neq -1.} \\
\end{array}%
\right.$$ Thus we obtain \be\label{atpalo}\Delta A_j = \left\{%
\begin{array}{ll}
\overline{Q}\bmatrix{ 0 & (\frac{(\overline{\lam})^{m-j}}{\|\Lam\|_2^2} Q_1^T r)^{T} \\
  \frac{(\overline{\lam}^j}{\|\Lam\|_2^2} Q_1^T r & X_j}Q^{H} \hfill{\, \mbox{if}\, \lam=-1.} \\
\overline{Q}\bmatrix{\frac{(\overline{\lam})^{j}+(\overline{\lam})^{m-j}}{2\|\Pi_+(\Lam)\|_2^2}x^Tr & (\frac{(\overline{\lam})^{m-j}}{\|\Lam\|_2^2} Q_1^T r)^{T} \\
\frac{(\overline{\lam})^j}{\|\Lam\|_2^2} Q_1^T r & X_j }Q^{H}, \hfill{\, \mbox{if} \, \lam\neq -1.}\end{array}%
\right.\ee Now setting $X_j=0, j=0:m$ and by (\ref{qtr}) we obtain $$\eta_F^{\S_p}(\lam,x,\Pb)=  \left\{%
\begin{array}{ll}
\frac{\|r\|_2}{\sqrt{m+1}}, \hfill{\, \mbox{if} \, \lam= -1} \\
\sqrt{\frac{2}{\|\Lam\|_2^2} \|r\|_2^2 + \frac{\|\Pi_-(\lam)\|_2^2 -\|\Pi_+(\lam)\|_2^2 }{\|\Lam\|_2^2 \|\Pi_+(\lam)\|_2^2}|x^Tr|^2}, \hfill{\, \mbox{if} \, \lam\neq -1}
\end{array}%
\right.$$ Moreover by DKW Theorem, (\ref{atpalo}) and following a similar techenique used for even $m,$ we obtain $$\eta_2^{\S_p}(\lam,x,\Pb) = \left\{%
\begin{array}{ll}
\sqrt{\frac{4\|\Pi_+(\Lam) + \Pi_-(\Lam)\|_2^2 }{\|\Lam\|_2^4}\|r\|_2^2 + \frac{\|\Lam\|_2^4 -4\|\Pi_+(\Lam) + \Pi_-(\Lam)\|_2^2 \|\Pi_+(\Lam)\|_2^2}{\|\Lam\|_2^4 \|\Pi_+(\Lam)\|_2^2}|x^Tr|^2},\\ \hfill{\, \mbox{if} \, |\lam|> 1} \\
\sqrt{\frac{4\|\Pi_+(\Lam) - \Pi_-(\Lam)\|_2^2 }{\|\Lam\|_2^4}\|r\|_2^2 + \frac{\|\Lam\|_2^4 -4\|\Pi_+(\Lam) - \Pi_-(\Lam)\|_2^2 \|\Pi_+(\Lam)\|_2^2}{\|\Lam\|_2^4 \|\Pi_+(\Lam)\|_2^2}|x^Tr|^2},\\ \hfill{\, \mbox{if} \, |\lam|\leq 1.} \\
\end{array}%
\right.$$ Hence the result follows when $\S=\S_p$ and $m$ is odd. Following a similar arguments we obtain the desired results for $\S=\S_{ap}.$ $\blacksquare$

\begin{remark}
Observe from the above proof that $\eta_F^\S(\lam,x,\Pb)$ is obtained by the only choice $X_j=0.$ For $\eta_2^\S(\lam,x,\Pb),$ by DKW Theorem, the choice of $X_j$ is infinite. Therefore the minimal structured perturbation is unique for Frobenious norm and in contrast we have infinitely many minimal structured perturbations for spectral norm.
\end{remark}
Let $\Pb\in\S.$ Treating $(\lam,x)\in\C\times\C^n$ with $\|x\|_2=1$ as an approximate eigenpair of $\Pb,$ we now construct a minimal structured perturbation $\Delta\Pb$ by simplifying the expressions of $\Delta A_j$ given in the proof of Theorem \ref{bck:tpal}. Let $P_x := I-xx^H$ where $I$ is the identity matrix of order $n$ and $0\neq x\in\C^n.$ Define \beano E_j &:=& \frac{(\overline{\lam})^{j}+\ep(\overline{\lam})^{m-j}}{2\|\Pi_{s}(\Lam)\|_2^2}(x^Tr)\overline{x}x^H, \,\,\, G_j := \frac{(\overline{\lam})^j +\ep(\overline{\lam})^{m-j}}{2\|\Pi_s(\Lam)\|_2^2 -\alpha |\lam^{m/2}|^2}(x^Tr)\overline{x}x^H + F_j \\ F_j &:=& \frac{1}{\|\Lam\|_2^2}\big[(\overline{\lam})^j P_x^T rx^H+ \ep(\overline{\lam})^{m-j} \overline{x}r^T P_x\big], \,\,\,H_{m/2} := \frac{(\overline{\lam})^{m/2} \, (x^Tr)\overline{x}x^H}{2\|\Pi_s(\Lam)\|_2^2 -\alpha |\lam^{m/2}|^2} + F_{m/2} \\ K_j &:=& \frac{(|\lam^j|^2(\overline{\lam})^{m-j}+\ep|\lam^{m-j}|^2(\overline{\lam})^j)~\overline{x^Tr}P_x^T rr^T P_x}{\|r\|_2^2-|x^Tr|^2},\,\,  L_j:=\frac{(\overline{\lam})^{m-j}~\overline{x^Tr}P_x^T rr^T P_x}{\|\Lam\|_2^2(\|r\|_2^2-|x^Tr|^2)}\eeano where $s\in\{ +, -\}, \ep\in\{+1, -1\}, \alpha \in \{0,1\}$ and $j\in\{0,1, \ldots, m\}.$ 
\begin{corollary}\label{cor:tpal}
Let $\S\in\{\S_p, \S_{ap}\}$ and $\Pb\in\S.$ Let $(\lam,x)$ be an approximate eigenpair of $\Pb.$ Then the unique structured perturbation $\Delta\Pb\in\S$ when $M=F,$ and a structured perturbation $\Delta\Pb\in\S$ when $M=2,$ of $\Pb$ for which $\Pb(\lam)x + \Delta\Pb(\lam)x =0$ and $\nrm{\Delta\Pb}_M = \eta^\S_M(\lam,x,\Pb)$ are given by 
%\begin{table}[h] 
\begin{enumerate} \item $m$ is odd:
\begin{center}\renewcommand{\arraystretch}{0.9}
%\begin{table}
\begin{tabular}{|l|c|c|}
%\multicolumn{3}{l}{$1. \,\,\, m$ is odd:} \\\\\hline
\hline
 & $\nrm{\cdot} = \nrm{\cdot}_F$ & $\nrm{\cdot} = \nrm{\cdot}_2$ \\
\hline
$\Delta A_j$ & $F_j$ if  $\lam=-1, \S=\S_p$ & $F_j$ if  $\lam=-1, \S=\S_p$ \\ \cline{2-3}
 & $E_j+F_j$ if $\lam\neq -1, \S=\S_p$ or & $E_j+F_j-\frac{|\lam^{m-j}|^{-2}}{2 \|\Pi_s(\Lam)\|_2^2}K_j$ if $|\lam| > 1$ \\
 & $\lam\neq 1, \S=\S_{ap}$  & and $\S\in\{\S_p, \S_{ap}\}$ \\ \cline{2-3}
& & $E_j+F_j-\frac{|\lam^{j}|^{-2}}{2 \|\Pi_s(\Lam)\|_2^2}K_j$ if $|\lam| \leq 1$ and \\
& $\frac{1}{\|\Lam\|_2^2}[rx^H - \overline{x}r^T]$ &  $\lam\neq -1, \S=\S_p$ or $\lam\neq 1, \S=\S_{ap}$  \\ \cline{3-3}
&  if $ \lam = 1, \S=\S_{ap}$ & $\frac{1}{\|\Lam\|_2^2}[rx^H - \overline{x}r^T]$ if $ \lam = 1, \S=\S_{ap}$ \\ 
\hline
\end{tabular}
%\end{table}\caption{\label{table1:tpalbckerr} $\Delta\Pb$ when $m$ is odd. }
\end{center}

\item $m$ is even:
\begin{center}\renewcommand{\arraystretch}{1.2}
\begin{tabular}{|l|c|c|}
% \multicolumn{3}{l}{$2. \,\,\, m$ is even:} \\\\\hline
\hline
 & $\nrm{\cdot} = \nrm{\cdot}_F$ & $\nrm{\cdot} = \nrm{\cdot}_2$ \\
\hline
$\Delta A_j$ & $\frac{\lam^j(x^Tr)}{\|\Lam\|_2^2}\overline{x}x^H + F_j$ if 
 & $\frac{\lam^j(x^Tr)}{\|\Lam\|_2^2}\overline{x}x^H + F_j-L_j$ if  \\ 
& $\lam = \pm 1, \S=\S_p$ & $\lam = \pm 1, \S=\S_p$ \\ \cline{2-3}
& $G_j$ if $\lam \neq \pm 1, \S=\S_p$  & $G_j-\frac{|\lam^{m-j}|^{-2}}{2 \|\Pi_s(\Lam)\|_2^2-\alpha |\lam^{m/2}|^2}K_j$  \\
& or $\lam \neq 1, \S=\S_{ap}$ & if $|\lam| > 1, \S\in\{\S_p,\S_{ap}\}$ \\ \cline{2-3}
& & $G_j-\frac{|\lam^{j}|^{-2}}{2 \|\Pi_s(\Lam)\|_2^2 - \alpha |\lam^{m/2}|^2}K_j$ if $|\lam|\leq 1$ and \\
& $\frac{1}{\|\Lam\|_2^2}[rx^H - \overline{x}r^T]$ &  $\lam\neq \pm 1, \S=\S_p$ or $\lam\neq 1, \S=\S_{ap}$ \\\cline{3-3}
&  if $ \lam = 1, \S=\S_{ap}$ & $\frac{1}{\|\Lam\|_2^2}[rx^H - \overline{x}r^T]$ if $ \lam = 1, \S=\S_{ap}$ \\
\hline \hline
$\Delta A_{m/2}$ & $\frac{(x^Tr)}{\|\Lam\|_2^2}\overline{x}x^H + F_{m/2}$  & $\frac{(x^Tr)}{\|\Lam\|_2^2}\overline{x}x^H + F_{m/2}-L_{m/2}$  \\
& if $\lam = \pm 1, \S=\S_p$ & if $\lam=\pm 1, \S=\S_p$\\\cline{2-3}
 & $H_{m/2}$ if $\lam \neq \pm 1, \S=\S_p$ & $H_{m/2}-\frac{\|\Lam\|_2^2}{2 \|\Pi_s(\Lam)\|_2^2-\alpha |\lam^{m/2}|^2}L_{m/2}$  \\
& & if $\lam\neq \pm 1, \S=\S_p$ \\\cline{2-3}
 & $\frac{(\overline{\lam})^{m/2}}{\|\Lam\|_2^2}[rx^H - \overline{x}r^T]$ if  & $\frac{(\overline{\lam})^{m/2}}{\|\Lam\|_2^2}[rx^H - \overline{x}r^T]$ if  \\
& $\lam \in\C, \S=\S_{ap}$ & $\lam \in\C, \S=\S_{ap}$ \\
\hline
\end{tabular}
%\caption{\label{table:tpalbckerr} }
\end{center}
%\end{table} 
\end{enumerate} where $\ep =\left\{%
\begin{array}{ll} 1 & \hbox{if $\S=\S_p$,} \\ -1 & \hbox{if $\S=\S_{ap}$,}\end{array}%
\right. s =\left\{%
\begin{array}{ll} + & \hbox{if $\S=\S_p$,} \\ - & \hbox{if $\S=\S_{ap}$,}\end{array}%
\right.$ $\alpha =\left\{%
\begin{array}{ll} 0 & \hbox{if $m$ is odd,} \\ 1 & \hbox{if $m$ is even.}\end{array}%
\right.$ and $\Delta A_{m-j}=\ep(\Delta A_j)^T.$
\end{corollary}
\noin\pf First consider $\S=\S_{p}.$ Let $m$ be even. If $\lam=\pm 1$ then simplifying (\ref{apm1f}) we have \beano \Delta A_j &=& \frac{\lam^j}{\|\Lam\|_2^2}(x^Tr)\overline{x}x^H + \frac{1}{\|\Lam\|_2^2}[\lam^j \overline{Q}_1Q_1^Trx^H + \lam^{m-j} \overline{x}r^TQ_1Q_1^H] + \overline{Q}_1X_jQ_1^H \\ &=& \frac{\lam^j}{\|\Lam\|_2^2}(x^Tr)\overline{x}x^H + \frac{1}{\|\Lam\|_2^2}[\lam^j P_x^Trx^H + \lam^{m-j} \overline{x}r^TP_x] + \overline{Q}_1X_jQ_1^H \\ &=& \frac{\lam^j}{\|\Lam\|_2^2}(x^Tr)\overline{x}x^H + F_j + \overline{Q}_1X_jQ_1^H, \eeano and if $\lam\neq \pm 1$ then simplifying (\ref{apmn1f}) and (\ref{apmn1mf}) we have \beano \Delta A_j &=& \frac{(\overline{\lam})^j +(\overline{\lam})^{m-j}}{2\|\Pi_s(\Lam)\|_2^2 - |\lam^{m/2}|^2}(x^Tr)\overline{x}x^H + F_j + \overline{Q}_1X_jQ_1^H \\ \Delta A_{m/2} &=& \frac{(\overline{\lam})^{m/2}}{2\|\Pi_s(\Lam)\|_2^2 - |\lam^{m/2}|^2}(x^Tr)\overline{x}x^H + F_{m/2} + \overline{Q}_1X_{m/2}Q_1^H. \eeano Now setting $X_j=0$ for $M=F$ we obtain the unique polynomial $\Delta\Pb\in\S$, and putting $X_j$ given in (\ref{x:pm1})-(\ref{x:npm11}) for $M=2$ we obtain $\Delta \Pb\in\S$ such that $\Pb(\lam)x + \Delta\Pb(\lam)x=0$ and $\nrm{\Delta\Pb}_M=\eta^\S_M(\lam,x,\Pb).$ The proof is similar when $m$ is odd, $\S=\S_p,$ and $\S=\S_{ap}, m$ is either even or odd. $\blacksquare$

%Next consider $\S=\S_{ap}.$ Then the desired result follows by simplifying

Note that if $Y \in \C^{n\times n}$  is such that $Yx =0$ and
$Y^Tx=0$ then $Y = (I-xx^H)^TZ(I-xx^H)$ for some matrix $Z.$ Hence
from the proof of Theorem~\ref{bck:tpal} and Corollary \ref{cor:tpal} we obtain that if
${\mathrm{K}}$ is a $T$-palindromic (resp. $T$-anti-palindromic) polynomial such that
$\Pb(\lam)x+ {\mathrm{K}}(\lam)x =0$ then ${\mathrm{K}}(z) =
\Delta \Pb(z) + (I-xx^H)^T {\mathrm{N}}(z) (I-xx^H)$ for some
$T$-palindromic (resp. $T$-anti-palindromic) matrix polynomial ${\mathrm{N}},$ where $\Delta \Pb$ is given in Corollary \ref{cor:tpal}.

\subsection{$H$-palindromic and $H$-anti-palindromic matrix polynomials}

We now consider the set of $H$-palindromic and $H$-anti-palindromic polynomials denoted by $\S_p$ and $\S_{ap},$ respectively. To derive the structured backward error of approximate eigenpair of a polynomial $\Pb\in\S,$ where $\S\in\{\S_p, \S_{ap}\},$ we proceed as follows. Let $z \in \C .$ Let us define the maps $\vec: \C \rightarrow \R^2 $ and $ \texttt{M} : \C \rightarrow \R^{2 \times 2} $ by \be\label{def:vec}\vec(z) = \bmatrix{\re (z) \\ \im (z)} \, \mbox{and} \, \texttt{M}(z) = \bmatrix{\re (z) & - \im (z) \\ \im (z) & \re (z)}.\ee Then we have the following Lemma.
\begin{lemma}
Let $z\in\C$ and $\Sigma = \bmatrix{1 & 0 \\ 0 & -1 }.$ Then the following hold.\\
$(i) \,\, \vec(\overline{z}) = \Sigma \, \vec (z). \,\, (ii)\,\, \vec(z_1 z_2) = \texttt{M}(z_1) \vec(z_2), \,z_1, z_2 \in \C.$ \,\, $(iii) \,\, \texttt{M}(\overline{z}) =
\texttt{M}(z)^T.$ 
\end{lemma}
\noin\pf The proof is obvious.$\blacksquare$

Note that $\Pb$ is $H$-palindromic polynomial if and only if $i\Pb, i :=\sqrt{-1}$ is $H$-anti-palindromic polynomial. Thus the map $H$-palindromic $\mapsto$ $H$-anti-palindromic is an isometric isomorphism. Also observe that $\eta_M^\S (\lam,x,\Pb) = \eta_M^{i\S} (\lam,x,i\Pb)$ where $i\S :=\{ i\Pb : \Pb \in \S\}.$ We denote the Moore-Penrose pseudoinverse of a matrix $A$ by $A^\dag.$

\begin{theorem}\label{bck:hpal}
Let $\Pb\in\S_p$ be given by $\Pb(z)=\sum_{j=0}^m z^jA_j.$ Let $(\lam,x)$ with $\|x\|_2=1$ be an approximate eigenpair of $\Pb.$ Set $r :=-\Pb(\lam)x$ and $\Lam = [1, \, \lam, \, \ldots, \, \lam^m]^T.$ Then we have \beano \eta_F^{\S_p}(\lambda,x,\Pb) &=& \left\{%
\begin{array}{ll}
    \frac{1}{\sqrt{m+1}}~\sqrt{2\|r\|_{2}^{2}-|r^{H}x|^{2}}\leq
\sqrt{2} \eta(\lam,x,\Pb), & \hbox{if $|\lam| =1$} \\
    \sqrt{(2\|\widehat{r}\|_2^2 - \ep|e^T_{(m/2)+1}\widehat{r}|^2) + 2\frac{\|r\|_2^2 -
|x^Hr|^2}{\|\Lam\|_2^2} }, & \hbox{if $|\lam| \neq 1.$} \\
\end{array}%
\right. \\ \eta_2^{\S_p}(\lambda,x,\Pb) &=&
\left\{%
\begin{array}{ll}
     \eta(\lam,x,\Pb) \hfill{\mbox{if} \, |\lam| =1}  \\
     \sqrt{(2\|\widehat{r}\|_2^2 - \ep|e^T_{(m/2) +1}\widehat{r}|^2) + \frac{(4\|\Pi_{+}(\Lam) +\Pi_-(\Lam)\|_2^2 - \ep|\lam^{m/2}|^2)
     (\|r\|_2^2-|x^Hr|^2)}
    {\|\Lam\|_2^4} }, \\ \hfill{\mbox{if} \, |\lam|>1} \\
    \sqrt{(2\|\widehat{r}\|_2^2 - \ep|e^T_{(m/2) + 1}\widehat{r}|^2) + \frac{(4\|\Pi_{+}(\lam) -\Pi_-(\Lam)\|_2^2 - \ep|\lam^{m/2}|^2)
    (\|r\|_2^2-|x^Hr|^2)}
    {\|\Lam\|_2^4} }, \\ \hfill{\mbox{if} |\lam|<1.} \\
\end{array}%
\right. \eeano where $\ep =\left\{%
\begin{array}{ll} 0 & \hbox{if $m$ is odd,} \\ 1 & \hbox{if $m$ is even,}\end{array}%
\right. \widehat{r}= \left\{%
\begin{array}{ll} \bmatrix{ H_0 & H_1 & \hdots & H_{(m-1)/2}}^\dag
\vec(x^Hr) \hfill{\,\mbox{if} \, m \, \mbox{is odd},} \\ \bmatrix{ H_0 & H_1 & \hdots & H_{(m-2)/2} & H_{m/2}}^\dag \vec(x^Hr) \hfill{\,\mbox{if} \, m \, \mbox{is even},}\end{array}%
\right.$ $$H_j = \bmatrix{
  \re \, (\lam^j) + \re \, (\lam^{m-j}) & -\im \, (\lam^j) + \im \, (\lam^{m-j})  \\
  \im \, (\lam^j) + \im (\lam^{m-j}) & \re \,( \lam^j) - \re \,(\lam^{m-j})}, \, j= \left\{%
\begin{array}{ll} 0:(m-1)/2 & \hbox{if $m$ is odd,} \\ 0:(m-2)/2 & \hbox{if $m$ is even,}\end{array}%
\right.$$ $H_{m/2} = \bmatrix{ \re (\lam^{m/2}) \\ \im  (\lam^{m/2})}$ whenever $m$ is even, and $e_j$ is the $j$-th column of the identity matrix.
\end{theorem}
\noin\pf First suppose that $m$ is even. By Theorem~\ref{exist:tpal} it is evident that there exists a polynomial $\Delta \Pb \in\S_p$ for which $\Delta \Pb(\lam)x +\Pb(\lam)x=0.$ Let $Q=[x, \, Q_1]$ be a unitary matrix where $x$ is given and $Q_1\in\C^{n\times (n-1)}$ is an isometry such that $Q_1^Hx=0.$ Define
\be\label{dahpal}\Delta A_j := Q\bmatrix{a_{jj} & a_j^H \\
  b_j & X_j} Q^H,~ \Delta A_{m/2} := Q\bmatrix{a_{(m/2)(m/2)} & a_{m/2}^H \\
  a_{m/2} & X_{m/2}}Q^H, ~j=0:(m-2)/2\ee and $\Delta A_{m-j} = (\Delta A_j)^H.$ Since $\Delta \Pb(\lam)x
+\Pb (\lam)x=0,$ we have
$$\bmatrix{\sum_{j=0}^{m} \lam^j a_{jj}\\
  \sum_{j=0}^{(m-2)/2} \lam^j b_j +\sum_{j=0}^{m/2} \lam^{m-j} a_j}=\bmatrix{x^Hr \\ Q_1^Hr }.$$ The minimum norm solution of $\sum_{j=0}^{(m-2)/2}
\lam^j b_j +\sum_{j=0}^{m/2} \lam^{m-j} a_j = Q_1^Hr$ is given by $b_j=\frac{(\overline{\lam})^j}{\|\Lam\|_2^2},\, a_j=\frac{(\overline{\lam})^{m-j}}{\|\Lam\|_2^2}.$ Note that for $|\lam|=1,$ we have $\overline{x^{H}r}=(\overline{\lambda})^{m} x^{H}r.$ Hence the
minimum norm solution of $\sum_{j=0}^{m} \lam^j a_{jj}=x^Hr$ is given by
$a_{jj}=\frac{(\overline{\lam})^j}{\|\Lam\|_2^2}x^Hr,~j=0:m.$ Therefore we have \be\label{a:hp1f}\Delta A_j = Q\bmatrix{\frac{(\overline{\lam})^j}{\|\Lam\|_2^2}x^{H}r & \frac{\lam^{m-j}}{\|\Lam\|_2^2}(Q_{1}^{H}r)^H\\ \frac{(\overline{\lam})^j}{\|\Lam\|_2^2}Q_{1}^{H}r &  X_j }Q^{H}, \Delta A_{m/2} = Q\bmatrix{\frac{(\overline{\lam})^{m/2}}{\|\Lam\|_2^2}x^{H}r & \frac{\lam^{m/2}}{\|\Lam\|_2^2}(Q_{1}^{H}r)^H\\ \frac{(\overline{\lam})^{m/2} }{\|\Lam\|_2^2}Q_{1}^{H}r &  X_{m/2}}Q^{H},\ee
$\Delta A_{m-j} = (\Delta A_j)^j, ~j=0:(m-2)/2,$ which gives,
$\eta_F^{\S_p}(\lambda,x,\Pb) =
\frac{1}{\sqrt{m+1}}~\sqrt{2\|r\|_{2}^{2}-|r^{H}x|^{2}}.$

If $|\lam| \neq 1,$ by Lemma \ref{app:1}, \ref{it:7}, the minimum norm solution of $\sum_{j=0}^{m} \lam^j a_{jj}=x^Hr$ is given by $a_{jj}= e_{j+1}^T \widehat{r}, \, j= 0: m/2$ where $\widehat{r} = \bmatrix{
H_0 & H_1 & \hdots & H_{(m-2)/2} & H_{m/2}}^\dag \vec(x^Hr) $ and $$ H_j = \bmatrix{
  \re \, (\lam^j) + \re \, (\lam^{m-j}) & -\im \, (\lam^j) + \im \, (\lam^{m-j})  \\
  \im \, (\lam^j) + \im (\lam^{m-j}) & \re \, (\lam^j) - \re \,(\lam^{m-j})}, j = 0:(m-2)/2,H_{m/2} = \bmatrix{\re (\lam^{m/2}) \\ \im (\lam^{m/2}) }.$$
Therefore we have \be\label{a:hpn1f}\Delta A_j = Q\bmatrix{e_{j+1}^T \widehat{r} & \lam^{m-j}\frac{(Q_{1}^{H}r)^H}{\|\Lam\|_2^2}\\ \frac{(\overline{\lam})^j Q_{1}^{H}r}{\|\Lam\|_2^2} &  X_j }Q^{H},
\Delta A_{m/2} = Q\bmatrix{ e_{(m/2)+1}^T \widehat{r} & \lam^{m/2}\frac{(Q_{1}^{H}r)^H}{\|\Lam\|_2^2} \\ \frac{(\overline{\lam})^{m/2} Q_{1}^{H}r}{\|\Lam\|_2^2} &  X_{m/2}}Q^{H}\ee and $\Delta A_{m-j} = (\Delta A_j)^H,~j=0:(m-2)/2.$ Setting $X_j =
0, j=0:m/2$ we obtain
$$\eta_F^{\S_p} (\lam,x,\Pb) = \sqrt{(2\|\widehat{r}\|_2^2 - |e_{(m/2)+1}^T
\widehat{r}|^2) + 2\frac{\|r\|_2^2 - |x^Hr|^2}{\|\Lam\|_2^2} }.$$

Next we derive the result for spectral norm. Note that for $|\lam|=1,$ by (\ref{a:hp1f}) and Theorem~\ref{dkw} we have $ \mu_{\Delta A_j}= \frac{\|r\|_2}{\|\Lam\|_2^2}$ and
\be\label{x:hpl1} X_j=-\frac{\lam^{m-j}~r^Hx~Q_1^Hr(Q_1^Hr)^H}{\|\Lam\|_2^2~(\|r\|_2^2-|r^Hx|^2)},~X_{m/2}=-\frac{\lam^{m/2}~r^Hx~Q_1^Hr(Q_1^Hr)^H}{\|\Lam\|_2^2~
(\|r\|_2^2-|r^Hx|^2)}.\ee Thus we have $ \eta_2^{\S_p}(\lambda,x,\Pb) =
\frac{\|r\|_{2}}{\sqrt{m+1}}.$ If $|\lam| \neq 1$ by (\ref{a:hpn1f}) and Theorem~\ref{dkw} we have $$\mu_{\Delta A_j} = \left\{%
\begin{array}{ll}
    \sqrt{|e_{j+1}^T\widehat{r}|^2+\frac{|\lam^{m-j}|^2~(\|r\|_2^2-|x^Hr|^2)}{\|\Lam\|_2^4}}, & \hbox{if $|\lam|>1,$} \\
    \sqrt{|e_{j+1}^T \widehat{r}|^2+\frac{|\lam^{i}|^2~(\|r\|_2^2-|x^Hr|^2)}{\|\Lam\|_2^4}}, & \hbox{if $|\lam|<1,$} \\
\end{array}%
\right. $$ $\mu_{\Delta A_{m/2}} = \sqrt{|e_{(m/2)+1}^T \widehat{r}|^2+\frac{|\lam^{m}|^2~(\|r\|_2^2-|x^Hr|^2)}{\|\Lam\|_2^4}}$  given by \be\label{x:hp1} X_j = \left\{%
\begin{array}{ll}
-\frac{\overline{e_{j+1}^T
\widehat{r}}~(\overline{\lam})^j~\lam^{m-j}~Q_1^Hr(Q_1^Hr)^H}{|\lam^{m-j}|^2~(\|r\|_2^2-|x^Hr|^2)},
& \hbox{if $|\lam|>1$} \\
    -\frac{\overline{e_{j+1}^T
\widehat{r}}~(\overline{\lam})^j~\lam^{m-j}~Q_1^Hr(Q_1^Hr)^H}{|\lam^{i}|^2~(\|r\|_2^2-|x^Hr|^2)}, & \hbox{if $|\lam|<1,$} \\
\end{array}%
\right.\ee \be\label{x:hp2} X_{m/2} =- \frac{(e_{(m/2)+1}^T
\widehat{r})~Q_1^Hr(Q_1^Hr)^H}{~(\|r\|_2^2-|x^Hr|^2)}\ee for
$j=0:(m-2)/2.$ This gives $$\eta_2^{\S_p}(\lam,x,\Pb) = \left\{%
\begin{array}{ll}
    \sqrt{(2 \|\widehat{r}\|_2^2 - |e_{(m/2)+1}^T
\widehat{r}|^2) + \frac{(4\|\Pi_{+}(\Lam) +\Pi_-(\Lam)\|_2^2 -
|\lam^{m/2}|^2) (\|r\|_2^2
- |x^Hr|^2)}{\|\Lam\|_2^4}  }, \\ \hfill{\mbox{if} |\lam|>1} \\
    \sqrt{(2 \|\widehat{r}\|_2^2 - |e_{(m/2)+1}^T
\widehat{r}|^2) + \frac{(4\|\Pi_{+}(\Lam)-\Pi_-(\Lam)\|_2^2 - |\lam^{m/2}|^2)
(\|r\|_2^2- |x^Hr|^2)}{\|\Lam\|_2^4}  }, \\ \hfill{\mbox{if} |\lam|<1.} \\
\end{array}%
\right.$$ Note that if $|x^Hr|=\|r\|_2,$ then $\|Q_1^Hr\|_2=0$. In such a case,
considering $X_j=0$ we obtain the desired results. Hence we are done for even $m.$ Following similar arguments the desired result can be obtained whenever $m$ is odd.$\blacksquare$

\begin{remark}
Observe from the above proof that $\eta_F^{\S_p}(\lam,x,\Pb)$ is obtained by the only choice $X_j=0.$ For $\eta_2^{\S_p}(\lam,x,\Pb),$ by DKW Theorem, the choice of $X_j$ is infinite. Therefore the minimal structured perturbation is unique for Frobenious norm and in contrast we have infinitely many minimal structured perturbations for spectral norm.
\end{remark}

Let $(\lam,x)$ with $\|x\|_2=1$ be an approximate eigenpair of a polynomial $\Pb\in\S_p.$ Now we construct a minimal structured perturbation $\Delta\Pb$ by simplifying $\Delta A_j$ given in the proof of Theorem \ref{bck:hpal}. We proceed as follows. Let $\tilde{E}_j := \frac{(\overline{\lam})^j}{\|\Lam\|_2^2}(r^{H}x)xx^H, \, \tilde{F}_j := \frac{1}{\|\Lam\|_2^2}[\lam^{m-j} xr^H P_x +(\overline{\lam})^j  P_x rx^H]$ and $  \tilde{K} = \frac{P_xrr^HP_x}{\|r\|_2^2 - |x^Hr|^2}.$ 

\begin{corollary}\label{cor:hpal}
Let $\Pb\in\S_p$ and $(\lam,x)$ be an approximate eigenpair of $\Pb.$ Then the unique structured perturbation $\Delta\Pb\in\S$ when $M=F,$ and a structured perturbation $\Delta\Pb$ when $M=2,$ of $\Pb$ for which $\Pb(\lam)x + \Delta\Pb(\lam)x =0$ and $\nrm{\Delta\Pb}_M = \eta^{\S_p}_M(\lam,x,\Pb)$ are given by 
%\begin{table}[h]
\begin{center}\renewcommand{\arraystretch}{0.9}
\begin{tabular}{l|c|c}
 & $\nrm{\cdot} = \nrm{\cdot}_F$ & $\nrm{\cdot} = \nrm{\cdot}_2$ \\
\hline
$\Delta A_j$ & $\tilde{E}_j + \tilde{F}_j$ if  $|\lam|=1$ & $\tilde{E}_j + \tilde{F}_j - x^Hr\lam^{m-j}\|\Lam\|_2^{-2} \tilde{K}$ if  $|\lam|=1$ \\
 & $e_{j+1}^T \widehat{r} xx^H + \tilde{F}_j$ if $|\lam|\neq 1$ & $e_{j+1}^T \widehat{r} xx^H + \tilde{F}_j - \overline{e_{j+1}^T \widehat{r}} (\overline{\lam})^j \lam^{m-j} |\lam^{m-j}|^{-2} \tilde{K}$ if $|\lam| > 1$ \\
& & $e_{j+1}^T \widehat{r} xx^H + \tilde{F}_j - \overline{e_{j+1}^T \widehat{r}} (\overline{\lam})^j \lam^{m-j} |\lam^{j}|^{-2}\tilde{K}$ if $|\lam| < 1$ \\
\hline
\end{tabular}
%\caption{\label{table:tpalbckerr} }
\end{center}
%\end{table} 
where $j=0:(m-1)/2$ if $m$ is odd, and $j=0:m/2$ if $m$ is even.
\end{corollary}
\noin\pf Setting $X_j=0$ when $M=F$ and putting $X_j$ given in (\ref{x:hpl1})-(\ref{x:hp2}) when $M=2$ the proof follows by simplifying (\ref{a:hp1f}).$\blacksquare$

It is evident from the proof of Theorem~\ref{bck:hpal} and by Corollary \ref{cor:hpal} that,
if ${\mathrm{K}}$ is a $H$-palindromic matrix polynomial such that
$\Pb(\lam)x+ {\mathrm{K}}(\lam)x =0$ then ${\mathrm{K}}(z) =
\Delta \Pb(z) + (I-xx^H) {\mathrm{N}}(z) (I-xx^H)$ for some
$H$-palindromic polynomial ${\mathrm{N}},$ where $\Delta \Pb$ is
given in Corollary~\ref{cor:hpal}.

\section{Structured backward error and palindromic linearizations}\label{sec3}

The classical way to solve polynomial eigenvalue problem is to convert the polynomial $\Pb$ into an equivalent linear polynomial $\Lb,$ called a linearization of $\Pb,$ and compute the eigenelements of $\Lb.$ For a polynomial $\Pb$ the set of linearizations form a vector space $\L_1(\Pb)$ defined by \cite{mackey:thesis} \be \L_1(\Pb):= \{\Lb(\lam) : \Lb(\lam)(\Ls\otimes I) = v\otimes \Pb(\lam)\}, v\in\C^m, \ee where $v$ is called the right ansatz vector, $\otimes$ is the Kronecker product, $\Ls:=[\lam^{m-1}, \, \hdots \, \lam, \, 1]^T$ and $\Lb$ is of the form $\Lb(\lam)=\lam X + Y, X\in\C^{mn\times mn}, Y\in\C^{mn\times mn}.$ But an arbitrary linearizaton $\Lb\in\L_1(\Pb)$ can destroy the eigensymmetry of a structured polynomial $\Pb$ \cite{mackey:thesis}. Hence to solve a structured polynomial eigenvalue problem one needs to choose a linearization which preserves the eigensymmetry of the polynomial. These linearizations are called structured linearizations.

Mackey et al.\cite{mackey2} have shown that a $*$-palindromic/$*$-anti-palindromic matrix polynomial can have both $*$-palindromic and $*$-anti-palindromic linearizations that preserve the eigensymmetry of the polynomial. Table~\ref{table:palansatzvectors} gives the structure of ansatz
vectors for structured linearizations, where $R = \left[%
\begin{array}{ccc}
  &  & 1 \\
   &  {\mathinner{\mkern2mu\raise1pt\hbox{.}\mkern2mu
\newline \raise4pt\hbox{.}\mkern2mu\raise7pt\hbox{.}\mkern1mu}} &    \\
  1 \\
\end{array}%
\right],$ for details see~\cite{mackey:thesis,mackey2}.
\begin{table}[h]
\begin{center}\renewcommand{\arraystretch}{1.2}
\begin{tabular}{|c|c|l|}
  \hline
  $\S$ & Structured Linearization & ansatz vector \\
  \hline
  \hline
  $T$-palindromic & $T$-palindromic & $Rv=v$ \\\cline{2-3}
      & $T$-anti-palindromic & $Rv=-v$ \\
   \hline
  $T$-anti-palindromic & $T$-palindromic & $R v = -v$ \\ \cline{2-3}
   & $T$-anti-palindromic & $R v =v$ \\
  \hline
  \hline
  $H$-palindromic & $H$-palindromic & $Rv=\overline{v}$ \\\cline{2-3}
      & $H$-anti-palindromic & $Rv=-\overline{v}$ \\
   \hline
  $H$-anti-palindromic & $H$-palindromic & $R v = -\overline{v}$ \\ \cline{2-3}
   & $H$-anti-palindromic & $R v =\overline{v}$\\
   \hline
\end{tabular}
\caption{\label{table:palansatzvectors} Table of the admissible
ansatz vectors for palindromic polynomials.}
\end{center}
\end{table} 

Note that to solve palindromic polynomial eigenvalue problem, the prime task is to detect potential structured linearizations that behave well during computations. Analyzing sensitivity of eigenvalues, potential structured linearizations have been produced in \cite{ba:poly} for  $T$-polynomial/$T$-anti-palindromic matrix polynomials. With a view to analyze accuracy of approximate eigenelements, in this section, we identify the potential structured linearizations of a $*$-polynomial/$*$-anti-palindromic polynomial. 

We first review some basic results available in the literature. Let $\Pb$ be a regular polynomial. Let $\Lb\in\L_1(\Pb)$ be a linearization of $\Pb$ corresponding to the right ansatz vector $v\in\C^{m}.$ Then the relationship between the eigenelements of $\Pb$ with that of its linearization $\Lb$ is given in \cite{higg1,higcon} \begin{itemize} \item $x\in
\C^n$ is a right eigenvector for $\Pb $ corresponding to an
eigenvalue $\lam \in \C$ if and only if $\Ls \otimes x$ is an
eigenvector for $\Lb(\lam)$ corresponding to the eigenvalue
$\lam.$ \end{itemize} Treating $(\lam,x)$ as an approximate eigenpair of $\Pb$ the relations \begin{eqnarray} \label{ln1} \|\Lb(\lam) (\Ls \otimes
x)\|_2 &=& \|v\|_2 \, \|\Pb(\lam)x\|_2, \\ \label{ln2}  
|(\Ls \otimes x)^T
\Lb(\lam) (\Ls \otimes x)| &=& |\Ls^T v| \, |x^T \Pb(\lam)x|, \\ \label{ln3} 
|(\Ls \otimes x)^H \Lb(\lam) (\Ls \otimes x)| &=& |\Ls^H v| \,
|x^H \Pb(\lam)x|\end{eqnarray} have been derived in \cite{higg1}. In view of (\ref{ln1})-(\ref{ln3}), without loss of generality we assume that the right ansatz vector $v$ is of unit norm.   Note the inequality \be\sqrt{\frac{m+1}{2m}} \leq
\frac{ \|\Lam\|_2
}{\|\Ls\|_2 \, \|(\lam, \, 1)\|_2} \leq 1\label{bckerr:insig}\ee given in \cite[Lemma A.1]{higcon}.

We denote the backward error of $(\lam, \Ls\otimes x)$ by $\eta(\lam,\Ls\otimes x,\Lb;v).$  Comparing $\eta_M(\lam,x,\Pb)$ with $\eta_M(\lam,\Ls\otimes x,\Lb;v)$ the inequality \be\label{unstr:bck} \sqrt{\frac{m+1}{2m}} \leq \frac{\eta_M (\lam,\Ls
\otimes x,\Lb;v)}{\eta_M (\lam,x,\Pb)} \leq 1 \ee has been proved in \cite[Theorem 4.1]{ba:poly}.  Now recall that $\eta_M(\lam,x,\Pb) \leq \eta_M^\S(\lam,x, \Pb).$ Hence for any structured linearization $\Lb\in\L_1(\Pb)$ of a given $\Pb\in\S,$ we have $\eta_M(\lam,\Ls\otimes x,\Lb;v) \leq \eta_M^\S(\lam,\Ls\otimes x, \Lb;v)$ where the
ansatz vector $v$ is given in Table~\ref{table:palansatzvectors}. Thus by \cite[Lemma 4.2]{ba:poly} we have \be\label{unbck:lemma} \sqrt{\frac{m+1}{2m}} \leq \frac{\eta_M^\S(\lam,\Ls\otimes x,\Lb;v)}{\eta_M(\lam,x,\Pb)}, M\in\{F,2\}.  \ee In the sequel we use the inequality \be\label{ln5} 0 < \frac{|\Ls^*v|}{\|\Ls\|_2}\leq 1, *\in\{T, H\}\ee where $v\in\C^m$ with $\|v\|_2=1.$

\begin{theorem}\label{comp:tpalin} Let $\Pb$ be a $T$-palindromic
 matrix polynomial. Let $\S_p\subset\L_1(\Pb)$ and $\S_{ap}\subset\L_1(\Pb)$ be the set of $T$-palindromic and $T$-anti-palindromic linearizations of $\Pb,$ respectively. Suppose $\Lb_{p} \in \S_p $ and $\Lb_{ap}\in\S_{ap}$ are the $T$-palindromic linearization and
$T$-anti-palindromic linearization of $\Pb$ corresponding to the
ansatz $Rv=v$ and $Rv=-v,$ respectively. If $(\lam,x)$ with $\|x\|_2=1$ is an approximate right eigenpair of $\Pb$ then we have
\begin{itemize} \item $M=F:$
\begin{enumerate} \small{ \item $\re(\lam)\geq 0:$ \,\,\,$\sqrt{1 -
\frac{2\, \re(\lam)}{|1+\lam|^2} } \sqrt{\frac{m+1}{m}}
\leq \dfrac{\eta^{\S_p}_F(\lambda,\Ls\otimes
x,\Lb_{p};v)}{\eta(\lam,x,\Pb)} \leq \sqrt{2}$ \item
$\re(\lam)\leq 0:$ \,\,\, $\sqrt{1 + \frac{2\, \re(\lam)}{|1-\lam|^2} } \sqrt{\frac{m+1}{m}} \leq
\dfrac{\eta^{\S_{ap}}_F(\lambda,x,\Lb_{ap};v)}{\eta(\lam,x,\Pb)} \leq
\sqrt{2}$.}
\end{enumerate}
\item $M=2:$ \begin{enumerate}\small{
\item $\lam\neq -1:$ \,\,\,
$\sqrt{\frac{m+1}{2m}} \leq \dfrac{\eta^{\S_{p}}_2(\lambda,x,\Lb_{p};v)}{\eta(\lam,x,\Pb)} \leq
\sqrt{2}~\sqrt{1+ \frac{1 }{|1+\lam|^2}}$
\item $\lam\neq 1:$ \,\,\,
$\sqrt{\frac{m+1}{2m}} \leq \dfrac{\eta^{\S_{ap}}_2(\lambda,x,\Lb_{ap};v)}{\eta(\lam,x,\Pb)}
\leq  \sqrt{2}~\sqrt{1+ \frac{1}{|1-\lam|^2}}.$}
\end{enumerate}
\end{itemize}
\end{theorem}

\noin\pf Let $r:=-\Pb(\lam)x.$ First consider $M=F.$ By Theorem \ref{bck:tpal} and using (\ref{ln1})-(\ref{ln2}) we have \beano \frac{\eta^{\S_p}_F(\lambda,x,\Lb_{p};v)}{\eta_M(\lam,x,\Pb)} &=&
\sqrt{2}~\frac{ \|\Lam\|_2 \sqrt{ \|r\|_2^2-2
\frac{\re(\lam)}{|1+\lam|^2} \,  \frac{|\Ls^Tv|^2}{\|\Ls\|_2^2} \,
|x^Tr|^2}}{ \|r\|_2 \|(1,\,\lam)\|_2 \, \|\Ls\|_2}\eeano if $\lam\neq -1,$ and $$\frac{\eta^{\S_{ap}}_F(\lambda,x,\Lb_{ap};v)}{\eta_M(\lam,x,\Pb)} =\sqrt{2}~\frac{\|\Lam\|_2 \sqrt{ \|r\|_2^2+ 2
\frac{\re(\lam)}{|1-\lam|^2} \,  \frac{|\Ls^Tv|^2}{\|\Ls\|_2^2} \,
|x^Tr|^2}}{\|r\|_2 \|(1,\,\lam)\|_2 \, \|\Ls\|_2} $$ if $\lam\neq 1.$ It is easy to verify that $\sqrt{1 - \frac{2 \re(\lam)}{|1 + \lam|^2}} \|r\|_2 \leq \sqrt{ \|r\|_2^2 - 2 \frac{\re(\lam)}{|1+\lam|^2} \,  \frac{|\Ls^Tv|^2}{\|\Ls\|_2^2} \,
|x^Tr|^2} \leq \|r\|_2$ if $\lam\neq -1,$ and  $\sqrt{1 + \frac{2 \re(\lam)}{|1 - \lam|^2}} \|r\|_2 \leq \sqrt{ \|r\|_2^2 + 2 \frac{\re(\lam)}{|1-\lam|^2} \,  \frac{|\Ls^Tv|^2}{\|\Ls\|_2^2} \,
|x^Tr|^2} \leq \|r\|_2$ if $\lam\neq 1.$ Thus by (\ref{ln5}),(\ref{bckerr:insig}) and (\ref{unbck:lemma}) the desired result follows for $M=F.$ 

Now consider $M=2.$ If $|\lam| > 1$ then by Theorem \ref{bck:tpal}, (\ref{ln1})-(\ref{ln2}) we have  \beano \frac{\eta^{\S_p}_2(\lambda,\Ls \otimes
x,\Lb_{p};v)}{\eta_M(\lam,x,\Pb)} &\leq&
    \frac{\sqrt{2} \, \|\Lam\|_2}{\|\Ls\|_2 \, \|(1, \, \lam)\|_2}~\sqrt{ \frac{|\Ls^Tv|^2 }{|1+\lam|^2\|\Ls\|_2^2} +
    \frac{|\lam|^2}{ \|(1, \, \lam)\|_2^2}} \\
   &\leq&
    \sqrt{2}~\sqrt{ 1+ \frac{1 }{|1+\lam|^2}}\eeano by (\ref{ln5}) and~(\ref{bckerr:insig}). Hence the results follow by (\ref{unbck:lemma}). Similarly if $|\lam| \leq 1, \lam\neq 1,$ by Theorem~\ref{bck:tpal} we have $$\frac{\eta^{\S_{ap}}_2(\lambda,\Ls\otimes
x,\Lb_{ap};v)}{\eta_M(\lam,x,\Pb)} \leq \sqrt{2}~\sqrt{ 1 + \frac{1}{|1-\lam|^2} }.$$ Hence the result follows by (\ref{unbck:lemma}). $\blacksquare$ 

Note that $|1-\lam|\leq |1+\lam|$ when $\re(\lam) \geq 0$ and $|1+\lam|\leq |1-\lam|$ when $\re(\lam) \leq 0.$  

\begin{remark} Let $\Pb$ be a $T$-anti-palindromic polynomial. Then we obtain similar bounds from Theorem \ref{comp:tpalin} by interchanging the role of $\Lb_p$ and $\Lb_{ap}.$\end{remark}

\begin{theorem}\label{com:strp}
Suppose that the assumptions of Theorem \ref{comp:tpalin} hold. Let
$\S$ be the set of $T$-palindromic polynomials. Then we have \begin{itemize} \item $M=F:$
\small{$\dfrac{\eta^{\S_p}_F(\lambda,\Ls\otimes
x,\Lb_{p};v)}{\eta^\S_F(\lam,x,\Pb)} \leq \sqrt{2}$ if $\re(\lam)\geq 0,$ $
\dfrac{\eta^{\S_{ap}}_F(\lambda,x,\Lb_{ap};v)}{\eta^\S_F (\lam,x,\Pb)} \leq
\sqrt{2}$ if $\re(\lam)\leq 0.$}
\item $M=2:$ \small{
$\dfrac{\eta^{\S_{p}}_2(\lambda,x,\Lb_{p};v)}{\eta^\S_F(\lam,x,\Pb)} \leq
\sqrt{2}$ if $\lam\neq -1,$
$\dfrac{\eta^{\S_{ap}}_2(\lambda,x,\Lb_{ap};v)}{\eta^\S_F(\lam,x,\Pb)}
\leq \sqrt{2}$ if $\lam\neq 1.$}
\end{itemize}
\end{theorem}
\noin\pf The proof is followed from the fact that $\eta_M(\lam,x,\Pb) \leq \eta^\S_M(\lam,x,\Pb)$ and Theorem \ref{comp:tpalin}.$\blacksquare$

\begin{remark} Let $\S$ be the set of $T$-anti-palindromic polynomials and $\Pb\in\S.$ Then the similar bounds hold for $T$-palindromic and $T$-anti-palindromic linearizations but the role of $T$-palindromic linearizations and $T$-anti-palindromic linearizations get exchanged in Corollary \ref{com:strp}.\end{remark}

The moral of the Theorem~\ref{comp:tpalin} and Corollary \ref{com:strp} is as follows. For a $T$-palindromic polynomial the bounds derived above advice to choose $T$-palindromic linearization
when $\re(\lam)\geq 0,$ and choose $T$-anti-palindromic linearization when $\re(\lam)\leq 0.$  Observe that our choice of structured linearizations is compatible with that given in~\cite{ba3} by analyzing structured condition number.

Now we consider $H$-palindromic/$H$-anti-palindromic matrix polynomials.

\begin{theorem}\label{comp:htpalin} Let $\S_p$ be the set of $H$-palindromic
matrix polynomials and $\Pb\in\S_p$. Let $\S\subset\L_1(\Pb)$ be the set of $H$-palindromic or $H$-anti-palindromic linearizations of $\Pb.$ Suppose $\Lb \in \S $ is an $H$-palindromic linearization or $H$-anti-palindromic linearization of $\Pb$ corresponding to the
ansatz $Rv=\overline{v}$ or $Rv=-\overline{v},$ respectively. If $(\lam,x)$ with $\|x\|_2 = 1$ is an approximate right eigenpair of $\Pb$ and $|\lam|=1$ then we have
\beano \sqrt{\frac{m+1}{2m}} \leq \frac{\eta_F^\S(\lam,\Ls\otimes x;\Lb,v)}{\eta_F(\lam,x,\Pb)} \leq \sqrt{2} \,\, \mbox{and}\,\, \frac{\eta_F^\S(\lam,\Ls\otimes x;\Lb,v)}{\eta_F^{\S_p}(\lam,x,\Pb)} \leq \sqrt{2}  \\ \sqrt{\frac{m+1}{2m}} \leq \frac{\eta_2^\S(\lam,\Ls\otimes x;\Lb,v)}{\eta_2(\lam,x,\Pb)} = \frac{\eta_2^\S(\lam,\Ls\otimes x;\Lb,v)}{\eta_2^{\S_p}(\lam,x,\Pb)} \leq \sqrt{2}. \eeano
\end{theorem}

\noin\pf Let $r:=-\Pb(\lam)x.$ If $|\lam|=1,$ by Theorem~\ref{bck:hpal} we have 
$$\frac{\eta_F^\S(\lam,\Ls\otimes x;\Lb,v)}{\eta(\lam,x,\Pb)} = \frac{\|\Lam\|_2 \sqrt{2\|r\|_2^2 - \frac{|\Ls^Hv|^2 \, |x^Hr|^2}{\|\Ls\|_2^2}}}{\|r\|_2 \|(1, \,
\lam)\|_2 \, \|\Ls\|_2} .$$ It is easy to verify that $\|r\|_2 \leq \sqrt{2\|r\|_2^2 - \frac{|\Ls^Hv|^2 \, |x^Hr|^2}{\|\Ls\|_2^2}} \leq \sqrt{2} \|r\|_2.$ Hence the desired result follows by
(\ref{unbck:lemma}). For the spectral norm we obtain the desired result by noting that $\eta^\S(\lam,x,\Pb) = \eta(\lam,x,\Pb)$ by Theorem~\ref{bck:hpal}.$\blacksquare$

This shows that there is almost no adverse effect of structured linearization of a $H$-palindromic matrix polynomial on the backward errors of approximate eigenelements when the approximate
eigenvalues are on the unit disk. On the other hand, if $|\lam| \neq 1$, a little calculation gives the following bounds:
\beano \frac{\eta_F^\S(\lam,\Ls\otimes x;\Lb,v)}{\eta(\lam,x,\Pb)} &\leq&
\sqrt{2}\, \sqrt{1 + \frac{\|\widehat{r}\|_2^2}{\|r\|_2^2}} \,\, \mbox{if} \,\, |\lam| \neq 1 \\
\frac{\eta_2^\S(\lam,\Ls\otimes x;\Lb,v)}{\eta(\lam,x,\Pb)} &\leq& \left\{%
\begin{array}{ll}
    \sqrt{2} \, \sqrt{
\frac{\|\widehat{r}\|_2^2}{\|r\|_2^2} + \frac{1}{ \|(1,
\lam^{-1})\|_2^2} }, & \hbox{ if $|\lam|>1$} \\\\
   \sqrt{2} \,
\sqrt{ \frac{\|\widehat{r}\|_2^2}{\|r\|_2^2} + \frac{1}{ \|(1,
\lam)\|_2^2} }, & \hbox{if $|\lam|<1$,} \\
\end{array}%
\right. \eeano where $\widehat{r}:= \widehat{r}_p =\bmatrix{ 1+\re(\lam) & \im(\lam) \\
\im(\lam) &
1-\re(\lam)}^\dag \bmatrix{ \re(\Ls^Hv \, x^H\Pb(\lam)x) \\
\im(\Ls^Hv \, x^H\Pb(\lam)x) }$ for $H$-palindromic linearization
and $\widehat{r}:=\widehat{r}_{ap}=\bmatrix{ 1-\re(\lam) & -\im(\lam) \\
-\im(\lam)
&1+\re(\lam)}^\dag \bmatrix{ \re(\Ls^Hv \, x^H\Pb(\lam)x) \\
\im(\Ls^Hv \, x^H\Pb(\lam)x) }$ for $H$-anti-palindromic
linearization.

Therefore the moral for $H$-palindromic/$H$-anti-palindromic matrix polynomials is as follows. If eigenvalues are on the unit disk then it does not matter whether we choose $H$-palindromic or $H$-anti-palindromic linearization. However, for eigenvalues not on the unit disk, it may be a good idea to solve both $H$-palindromic or $H$-anti-palindromic linearizations and pick an eigenpair $(\lam, x)$ from $H$-palindromic or $H$-anti-palindromic linearization according as $r_p \leq r_{ap}$ or $r_{ap} \leq r_p$.

\section{Conclusion}

We have derived computable expression of structured backward errors of approximate eigenpairs of $*$-palindromic/$*$-anti-palindromic matrix polynomials. We mention that these expressions have an important role to play in analyzing stability of structured preserving algorithms. Finally structured backward errors have been used to determine potential structured linearizations of a $*$-palindromic/$*$-anti-palindromic matrix polynomial.

\section{Acknowledgments} This work was done while the author was a PhD student at IIT Guwahati. The author thanks his supervisor Rafikul Alam for his stimulating comments and discussions that have significantly improved the quality of the results. The author would like to thank the reviewers for their comments that help improve the manuscript. Thanks are also due to Daniel Kressner for inspiring discussions on palindromic eigenvalue problems. Finally, thanks are also due to Harish K. Pillai for his comments and suggestions to improve the manuscript.

\appendix 
\section{Appendix}

The following Lemma summarizes few auxiliary results required to prove the Theorems in section \ref{sec2}.

\begin{lemma}\label{app:1}
Let $\lam\in\C$ and $\Lam = [1, \, \lam, \, \ldots, \, \lam^m]^T\in\C^{m+1}.$ Let $x,r\in\C^n, \ep\in\{+1, -1\}$ and $s\in\{+, -\}.$ Then the solution of $\sum_{j=0}^m a_{jj}\lam^j = x^*r, a_{jj}\in\C, *\in\{T, H\}$ that minimizes $\sum_{j=0}^m |a_{jj}|^2$ is given by \begin{enumerate}
\item\label{it:1} $a_{jj}=0$ if $x^*r=0.$
\item\label{it:2} $a_{jj}= \frac{(\overline{\lam})^j}{\|\Lam\|_2^2}x^*r.$  
\item\label{it:3} If $m$ is odd, $a_{jj}=\ep \, a_{(m-j)(m-j)}, j=0:(m-1)/2$ then $a_{jj}= \frac{(\overline{\lam})^j + \ep(\overline{\lam})^{m-j}}{2\|\Pi_s(\Lam)\|_2^2} x^*r,$ whenever $s=+$ if $\ep=1,$ and $s=-$ if $\ep=-1.$
\item\label{it:4} If $m$ is even, $a_{jj}=a_{(m-j)(m-j)}, j=0:m/2$ then $a_{jj}= \frac{(\overline{\lam})^j + (\overline{\lam})^{m-j}}{2\|\Pi_s(\Lam)\|_2^2-|\lam^{m/2}|^2} x^*r$ and $a_{(m/2)(m/2)} = \frac{(\overline{\lam})^{m/2}}{2\|\Pi_s(\Lam)\|_2^2-|\lam^{m/2}|^2}x^*r.$ 
\item\label{it:5} If $m$ is even, $a_{jj}=-a_{(m-j)(m-j)}, j=0:m/2$ then $a_{jj}= \frac{(\overline{\lam})^j - (\overline{\lam})^{m-j}}{2\|\Pi_-(\Lam)\|_2^2} x^*r$ and $a_{(m/2)(m/2)}=0.$
\item\label{it:6} If $m$ is odd, $\overline{a}_{jj} = \ep a_{m-j,m-j} \, j=0:(m-1)/2$ and $|\lam|\neq 1$ then $a_{jj} = e_{j+1}^T \widehat{r}$ where $$ \widehat{r} = \bmatrix{ H_0 & H_1 & \hdots & H_{(m-1)/2}}^\dag \vec(x^*r), \, H_j = \bmatrix{ \re \, (\lam^j) + \ep \re \, (\lam^{m-j}) & -\im \, (\lam^j) + \ep \im \, (\lam^{m-j})  \\  \im \, (\lam^j) + \ep \im (\lam^{m-j}) & \re \, (\lam^j) - \ep \re \, (\lam^{m-j}) },$$ $e_j$ is the $j$-th column of the identity matrix and $\vec$ is defined in (\ref{def:vec}).
\item\label{it:7} If $m$ is even, $\overline{a}_{jj} = \ep a_{m-j,m-j} \, j=0:m/2$ and $|\lam|\neq 1$ then $a_{jj} = e_{j+1}^T \widehat{r}$ where $$ \widehat{r} = \bmatrix{ H_0 & H_1 & \hdots & H_{(m-1)/2} & H_{m/2} }^\dag \vec(x^*r),$$ $H_{m/2} = \bmatrix{ \re \, (\lam^{m/2})  \\  \im \, (\lam^{m/2}) }$ if $\ep=1,$ $H_{m/2} = \bmatrix{0 & -1 \\ 1 & 0} \bmatrix{ \re \, (\lam^{m/2})  \\  \im \, (\lam^{m/2}) }$ if $\ep=-1,$ and  $H_j, j=0:(m-2)/2$ is same as that given in \ref{it:6}.
\end{enumerate} where $\Pi_s(\Lam)$ is defined in (\ref{eq:symmetricprojector}).
\end{lemma}
\noin\pf The proof of \ref{it:1} and \ref{it:2} are obvious. Now consider \ref{it:3}. Let $m$ be odd. Then $$\sum_{j=0}^m a_{jj}\lam^j = x^*r \Rightarrow \sum_{j=0}^{(m-1)/2} (\lam^{j} + \ep \lam^{m-j}) a_{jj} = x^*r.$$ Hence the result follows by \ref{it:2}. The proof is similar for \ref{it:4} and \ref{it:5}. To prove \ref{it:6} we proceed as follows. Let $\ep=1.$ Apply the map $\vec$ at both sides of $\sum_{j=0}^{m} \lam^i a_{jj} = x^*r.$ This yields
\be\label{eq1} \vec (\sum_{j=0}^{m} \lam^i a_{jj}) = \vec(x^*r) \,\, \Rightarrow \,\, \sum_{j=0}^{m} \texttt{M}(\lam^i) \vec(a_{jj}) = \vec (x^*r)\ee where $\vec$ and $\texttt{M}$ are defined in (\ref{def:vec}). Employing the condition $\overline{a}_{jj} = a_{m-j,m-j} \, j=0:(m-1)/2$ on (\ref{eq1}) we have $$\sum_{j=0}^{(m-1)/2} (\texttt{M} (\lam^j)  + \texttt{M} (\lam^{m-j}) \Sigma ) \vec
(a_{jj}) = \vec(x^*r)$$ where $\Sigma = \bmatrix{1 & 0 \\0 & -1}.$ Thus we have $a_{jj} = e_{j+1}^T \widehat{r}$ where $$ \widehat{r} = \bmatrix{ H_0 & H_1 & \hdots & H_{(m-1)/2}}^\dag \vec(x^*r), \, H_j = \bmatrix{ \re \, (\lam^j) + \re \, (\lam^{m-j}) & -\im \, (\lam^j) + \im \, (\lam^{m-j})  \\  \im \, (\lam^j) + \im (\lam^{m-j}) & \re \, (\lam^j) - \re \, (\lam^{m-j}) },$$ $e_j$ is the $j$-th column of the identity matrix. The proof is similar for $\ep=-1$ and \ref{it:7}.$\blacksquare$

\begin{lemma}\label{app:2}
Let $\lam\in\C, x_j\in\C^{n-1}, j=0:m.$ Then the solution of $\sum_{j=0}^m \lam^j x_j=y, y\in\C^{n-1}$ that minimizes $\sum_{j=0}^m \|x_j\|_2^2$ is given by $x_j=\frac{(\overline{\lam})^j}{\|\Lam\|_2^2} y$ where $\Lam = [1, \, \lam, \, \ldots, \, \lam^m]^T.$ 
\end{lemma}
\noin\pf The proof follows by using Moore-Penrose pseudoinverse of $\Lam$.$\blacksquare$

\end{document}